\input amstex\documentstyle{amsppt}  
\pagewidth{12.5cm}\pageheight{19cm}\magnification\magstep1
\topmatter
\title Study of a $\bold Z$-form of the coordinate ring of a reductive group\endtitle  
\author G. Lusztig\endauthor
\address{Department of Mathematics, M.I.T., Cambridge, MA 02139}\endaddress
\thanks{Supported in part by the National Science Foundation}\endthanks
\endtopmatter   
\document
\define\udBB{\un{\dot{\BB}}}
\define\uBB{\un{\BB}}
\define\hUU{\hat{\UU}}
\define\hOO{\hat{\OO}}
\define\dBB{\dot{\BB}}

\define\dUU{\dot{\UU}}

\define\ua{\un a}
\define\ub{\un b}
\define\uc{\un c}
\define\ui{\un i}

\define\ur{\un r}
\define\uy{\un y}

\define\hD{\hat D}

\define\hm{\hat m}
\define\hmu{\hat\mu}

\define\da{\dagger}

\define\Lie{\text{\rm Lie }}

\define\si{\sim}

\define\sqc{\sqcup}

\define\qua{\quad}

\define\ha{\hat a}
\define\hb{\hat b}
\define\hc{\hat c}

\define\lb{\linebreak}

\define\eSb{\endSb}

\define\bin{\binom}
\define\op{\oplus}

\define\part{\partial}

\define\ra{\rangle}
\define\n{\notin}
\define\iy{\infty}
\define\m{\mapsto}
\define\do{\dots}
\define\la{\langle}

\define\lra{\leftrightarrow}

\define\sub{\subset}    

\define\T{\times}
\define\ti{\tilde}
\define\nl{\newline}
\redefine\i{^{-1}}

\define\un{\underline}

\define\ot{\otimes}

\define\bcc{\thickfracwithdelims[]\thickness0}

\define\Hom{\text{\rm Hom}}
\define\End{\text{\rm End}}
\define\Aut{\text{\rm Aut}}

\define\di{\diamond}

\define\a{\alpha}

\redefine\c{\chi}
\define\g{\gamma}
\redefine\d{\delta}
\define\e{\epsilon}
\define\et{\eta}
\define\io{\iota}
\redefine\o{\omega}
\define\p{\pi}
\define\ph{\phi}
\define\ps{\psi}
\define\r{\rho}

\redefine\t{\tau}
\define\th{\theta}
\define\k{\kappa}
\redefine\l{\lambda}
\define\z{\zeta}
\define\x{\xi}

\redefine\G{\Gamma}
\redefine\D{\Delta}

\redefine\L{\Lambda}

\define\boc{\bold c}

\define\bof{\bold f}
\define\hh{\bold h}

\define\kk{\bold k}

\define\mm{\bold m}

\define\oo{\bold o}

\redefine\xx{\bold x}

\define\BB{\bold B}
\define\CC{\bold C}

\define\NN{\bold N}
\define\OO{\bold O}
\define\PP{\bold P}
\define\QQ{\bold Q}

\define\UU{\bold U}

\define\ZZ{\bold Z}

\define\ca{\Cal A}

\define\cm{\Cal M}

\define\co{\Cal O}

\define\car{\Cal R}

\define\ct{\Cal T}

\define\fA{\frak A}

\define\fC{\frak C}

\define\fF{\frak F}
\define\fG{\frak G}

\define\fT{\frak T}

\define\ta{\ti a}

\define\tc{\ti c}
\define\td{\ti d}
\define\te{\ti e}

\define\tih{\ti h}

\define\tn{\ti n}

\define\uZ{\un Z}

\define\cir{\circ}

\define\BO{Bo}
\define\CH{C1}
\define\CHH{C2}
\define\JO{Jo}
\define\KO{Ko}
\define\QG{L1}
\define\LU{L2}
\define\SO{So}

\head Introduction\endhead
In his famous paper \cite{\CH} Chevalley associated to any root datum of adjoint type
 and to any field $k$ a certain group (now 
known as a Chevalley group) which in the case where $k=\CC$ was the usual adjoint Lie group over $\CC$ and which 
in the case where $k$ is finite led to some new families of finite simple groups. 
Let $\OO'$
 be the coordinate ring of a connected semisimple group $G$ over $\CC$ attached to 
a fixed (semisimple) root datum. In a sequel \cite{\CHH} to \cite{\CH}, Chevalley defined 
a $\ZZ$-form of $\OO'$. Later, another construction of such a $\ZZ$-form was proposed by 
Kostant  \cite{\KO}. Kostant notes that $\OO'$ can be viewed as a "restricted"
dual of the universal enveloping algebra $\UU$ of $\Lie G$; he defines a $\ZZ$-form 
$\UU_\ZZ$ of $\UU$ ("the Kostant 
$\ZZ$-form") and then defines the $\ZZ$-form $\OO'_\ZZ$ as the set of all elements in $\OO'$ 
which take integral values on $\UU_\ZZ$. Then for any commutative ring $A$ with $1$
he defines 
$\OO'_A$ as $A\ot\OO'_\ZZ$; this is naturally a Hopf algebra over $A$. 
It follows that the set $G_A$ of $A$-algebra homomorphisms $\OO'_A@>>>A$ has a natural
group structure. Thus the root datum gives rise to a family of groups $G_A$, one for each
$A$ as above. 

Unlike Chevalley's
approach which was based on a choice of a faithful representation of $G$, Kostant's approach is direct (no choices
involved) and generalizes to the quantum case. 

In this paper we develop the theory of Chevalley groups following Kostant's approach. We shall prove that:

(I) {\it if $A$ is an algebraically closed field then $\OO'_A$ is the coordinate algebra of a connected semisimple algebraic group over $A$ corresponding to the given root datum.}
\nl
(We treat the reductive case at the same time.) 
Note that (I) was stated without proof in \cite{\KO}. 

In this paper we note that Kostant's definition can be reformulated by replacing $\UU$ 
by a "modified enveloping 
algebra". The theory is then developed using extensively the theory of canonical bases of 
such modified enveloping algebras (presented in \cite{\QG}), coming from quantum groups. (See the Notes in 
\cite{\QG} for references to original sources concerning canonical bases.)

We now present the content of this paper in more detail.

Let $A$ be a fixed commutative ring with $1$ with a given invertible element $v\in A$.

In \S1 we recall the definition and some properties of the modified (quantized) enveloping algebra $\dUU_A$ over
$A$ and its canonical basis $\dBB$. We also define a "completion" $\hUU_A$ of $\dUU_A$ which consists of formal
(possibly infinite) $A$-linear combinations of elements in $\dBB$. We show that the multiplication and
"comultiplication" of $\dUU_A$ extend naturally to $\hUU_A$.

In \S2 we introduce some invertible elements $s'_{i,e}$ of $\hUU_A$ (where $i$ corresponds to a simple reflection
and $e=\pm1$). We show that conjugation by $s'_{i,e}$ restricted to $\dUU_A$ is essentially the action of a
generator in the braid group action on $\dUU_A$ studied in \cite{\QG}. Thus $s'_{i,e}$ plays the same role as an
element considered in a similar context (with $A=\CC(v)$) by Soibelman \cite{\SO}. But while Soibelman's element 
is not explicit and its integrality properties are not clear, our element $s'_{i,e}$ is remarkably simple and has 
obvious integrality properties.

In \S3 we define following \cite{\QG, 29.5.2} the Hopf algebra $\OO_A$ (a quantum analogue
of $\OO'_A$ above.) We prove that the $A$-algebra $\OO_A$ is finitely generated (see 3.3) with an explicit set of
generators. In 3.7 we show that the $A$-algebra $\OO_A$ can be imbedded into the tensor product of two simpler
algebras. (For a closely related result in the case where $A=\QQ(v)$, see \cite{\JO, 9.3.13}.) In the case
where $v=1$ in $A$ these simpler algebras can be explicitly described in terms of the braid group action, see
3.13. We deduce that, if $A$ is an integral domain and $v=1$ in $A$, then $\OO_A$ is an
integral domain; see Theorem 3.15. (For a similar result in the case where $A=\QQ(v)$, see
\cite{\JO, 9.1.9}.)

In \S4 we assume that $v=1$ in $A$ and we introduce the group $G_A$ in analogy with
\cite{\KO}. In Theorem 4.11 we show that $\OO_A$ has a property like (I) above.

In \S5 we identify (assuming that $v=1$ in $A$) our $\OO_A$ with Kostant's $\OO'_A$.

\head Contents\endhead
1. The algebras $\dUU_A$, $\hUU_A$.

2. The elements $s'_{i,e},s''_{i,e}$ of $\hUU_A$.

3. The Hopf algebra $\OO_A$.

4. The group $G_A$.

5. From enveloping algebras to modified enveloping algebras.

\head 1. The algebras $\dUU_A$, $\hUU_A$\endhead
\subhead 1.1\endsubhead
In this section we recall the definition of the modified quantized enveloping algebra
$\dUU_A$ (over $A$) attached to a root datum and we recall the definition and some of the
properties of the canonical basis $\dBB$ of $\dUU_A$.
We also study a certain completion $\hUU_A$ of $\dUU_A$.

We fix a root datum as in \cite{\QG, 2.2}. This consists of two free abelian groups of finite type $Y,X$ with a 
given perfect pairing $\la,\ra:Y\T X\to\ZZ$ and a finite set $I$ with given imbeddings $I@>>>Y\qua(i\m i)$ and 
$I@>>>X\qua(i\m i')$ such that $\la i,i'\ra=2$ for all $i\in I$ and $\la i,j'\ra\in-\NN$ for all $i\ne j$ in $I$;
in addition, we are given a symmetric bilinear form $\ZZ[I]\T\ZZ[I]\to\ZZ$, $\nu,\nu'\m\nu\cdot\nu'$ such that 
$i\cdot i\in2\ZZ_{>0}$ for all $i\in I$ and $\la i,j'\ra=2i\cdot j/i\cdot i$ for all $i\ne j$ in $I$. We assume 
that the matrix $(i\cdot j)_{i,j\in I}$ is positive definite. 

Let $X^+=\{\l\in X;\la i,\l\ra\ge0\text{ for all }i\in I\}$. For $\l,\l'$ in $X$ we write $\l\ge\l'$ or $\l'\le\l$
if $\l-\l'\in\sum_{i\in I}\NN i'$. The image of $\nu\in\ZZ[I]$ under the homomorphism $\ZZ[I]@>>>X$ such that 
$i\m i'$ for $i\in I$, is denoted again by $\nu$.

Let $W$ be the (finite) subgroup of $\Aut(Y)$ generated by the involutions 
$s_i:y\m y-\la y,i'\ra i\qua (i\in I)$ or 
equivalently the subgroup of $\Aut(X)$ generated by the involutions 
$s_i:x-\la i,x\ra i'\qua (i\in I)$; these two 
subgroups may be identified by taking contragredients. For $i\ne j$ in $I$ let $n_{ij}=n_{ji}$ be the order of 
$s_is_j$ in $W$. Let $l:W@>>>\NN$ be the standard length function on $W$ with respect to $\{s_i;i\in I\}$. Let 
$w_0\in W$ be the unique element such that $l(w_0)$ is maximal. 

\subhead 1.2\endsubhead
Let $v$ be an indeterminate. Let $\ca=\ZZ[v,v\i]$. For $i\in I$ we set $v_i=v^{i\cdot i/2}$. 

We fix a commutative ring $A$ with $1$ with a given ring homomorphism $\ca@>>>A$ respecting $1$. For $\a\in\ca$ we
shall often denote the image of $\a$ under $\ca@>>>A$ again by $\a$. Whenever we write $A=\QQ(v)$ or $A=\ca$ we 
shall understand that $A$ is regarded as an $\ca$-algebra in an obvious way. Whenever we write $A=\QQ$ or $A=\ZZ$
we shall understand that $A$ is regarded as an $\ca$-algebra with $v=1$ in $A$.

Let $A^\cir$ be the group of invertible elements of the ring $A$.

An $A$-linear map $\ph:H@>>>H'$ where $H,H'$ are $A$-modules is said to be a {\it split injection} if there exists
an $A$-linear map $\ps:H'@>>>H$ such that $\ps\ph=1$. 

\subhead 1.3\endsubhead
Let $\bof$ be the associative $\QQ(v)$-algebra with $1$ defined as in \cite{\QG, 1.2.5} or equivalently by the 
generators $\th_i(i\in I)$ and the quantum Serre relations
$$\sum_{p,p'\in\NN;p+p'=1-\la i,j'\ra}(-1)^{p'}(\th_i^p/[p]^!_i)\th_j(\th_i^{p'}/[p']^!_i)$$
for $i\ne j$ in $I$, where $[p]^!_i=\prod_{s=1}^p(v_i^s-v_i^{-s})/(v_i-v_i^{-1})$ for $p\in\NN$. We have a direct
sum decomposition $\bof=\op_{\nu\in\ZZ[I]}\bof_\nu$ as a vector space, where $\bof_\nu$ is spanned by words in 
$\th_i$ in which the number of apparitions of $\th_i$ is the coefficient of $i$ in $\nu$, for all $i\in I$. For 
$i\in I$, $n\in\ZZ$ we set $\th_i^{(n)}=\th_i^n/[n!]_i\in\bof$ if $n\ge0$ and $\th_i^{(n)}=0$ if $n<0$. Let 
$\bof_\ca$ be the $\ca$-subalgebra with $1$ of $\bof$ generated by the elements $\th_i^{(n)}$ with 
$i\in I,n\in\ZZ$. We have $\bof_\ca=\op_{\nu\in\ZZ[I]}\bof_{\nu,\ca}$ where $\bof_{\nu,\ca}=\bof_\nu\cap\bof_\ca$.
Let $\bof_A=A\ot_\ca\bof_\ca$, an $A$-algebra. We have $\bof_A=\op_{\nu\in\ZZ[I]}\bof_{\nu,A}$ where 
$\bof_{\nu,A}=A\ot_\ca\bof_{\nu,\ca}$.

Let $\BB$ be the canonical basis of $\bof$ (see \cite{\QG, 14.4}). Note that $\BB$ is also an $\ca$-basis of 
$\bof_\ca$. If $b\in\BB$ we shall denote the element $1\ot_\ca b\in A\ot_\ca\bof_\ca=\bof_A$ again by $b$. Thus 
$\BB$ can be viewed also as an $A$-basis of $\bof_A$.

\subhead 1.4\endsubhead
For any $n\in\ZZ,t\in\NN$ and $i\in I$ we define $\bcc{n}{t}_i\in\ca$ as in \cite{\QG, 1.3.1, 1.1.2}.

Let $\dUU_A$ be the $A$-algebra generated by the symbols $x^+1_\z x'{}^-$, $x^-1_\z x'{}^+$ with 
$x\in\bof_{\nu,A}$, $x'\in\bof_{\nu',A}$ for various $\nu,\nu'$ and $\z\in X$; these symbols are subject to the
following relations:

$\bof_A@>>>\dUU_A$, $x\m x^{\pm}1_\z$ is $A$-linear for any $\z\in X$;

$\th_i^{(n)+}1_\z\th_j^{(m)-}=\th_j^{(m)-}1_\z\th_i^{(n)+}$ if $m,n\in\NN$, $\z\in X$, $i\ne j$;

$\th_i^{(n)\pm}1_{\mp\z}\th_i^{(m)\mp}=\sum_{t\in\NN}\bcc{m+n-\la i,\z\ra}{t}_i
\th_i^{(m-t)\mp}1_{\mp\z\pm(a+b-t)i'}\th_i^{(n-t)\pm}$ if $m,n\in\NN$, $\z\in X$, $i\in I$;

$x^{\pm}1_\z=1_{\z\pm\nu}x^{\pm}$ if $x\in\bof_{\nu,A}$, $\z\in X$;

$(x^{\pm}1_\z)(1_{\z'}x'{}^{\mp})=\d_{\z,\z'}x^{\pm}1_\z x'{}^{\mp}$ if $x,x'\in\bof_A$, $\z\in X$;

$(x^{\pm}1_\z)(1_{\z'}x'{}^{\pm})=\d_{\z,\z'}1_{\z\pm\nu}(xx')^{\pm}$, if
$x\in\bof_{\nu,A},x'\in\bof_A$, $\z\in X$.
\nl
If $x$ or $x'$ in $x^+1_\z x'{}^-$ or $x^-1_\z x'{}^+$ is $1$ we omit writing it.
\nl
(See \cite{\QG, 31.1.1, 31.1.3}.) Note that the two $A$-linear maps $\bof_A\ot_AA[X]\ot_A\bof_A@>>>\dUU_A$, 
$x\ot\l\ot x'\m x^+1_\l x'{}^-$ and $x\ot\l\ot x'\m x^-1_\l x'{}^+$ are isomorphisms of $A$-modules. (Here $A[X]$ is
the free $A$-module with basis indexed by $X$.) Hence we have canonically $\dUU_A=A\ot_\ca\dUU_\ca$ as 
$A$-algebras.

Now $\dUU_A$ does not have $1$ in general; instead it has a family of 
elements $1_\l (\l\in X)$ such that $1_\l1_{\l'}=\d_{\l,\l'}$ for any $\l,\l'$ in $X$ and such that 
$\dUU_A=\sum_{\l,\l'\in X}1_\l\dUU_A1_{\l'}$ (necessarily a direct sum). 

In the case where $A=\QQ(v)$ we write $\dUU$ instead of $\dUU_A$. Now the obvious map $\dUU_\ca@>>>\dUU$ 
induced by the inclusion $\ca\sub\QQ(v)$ is injective and identifies $\dUU_\ca$ with a $\ca$-subalgebra of $\dUU$.
For any $\l,\l',\l_1,\l_2,\l'_1,\l'_2$ in $X$ we define an $A$-linear map
$$\D_{\l,\l',\l_1,\l'_1,\l_2,\l'_2}:1_\l\dUU_A1_{\l'}@>>>(1_{\l_1}\dUU_A1_{\l'_1})\ot_A(1_{\l_2}\dUU_A1_{\l'_2})$$
as follows: if $\l=\l_1+\l_2$, $\l'=\l'_1+\l'_2$ it is the map obtained by applying $A\ot_\ca?$ to the 
$\ca$-linear map at the end of \cite{\QG, 23.2.3}; otherwise, it is $0$.
We define an algebra isomorphism $S$ from $\dUU_A$ to $\dUU_A$ with opposed multiplication as follows: if
$A=\QQ(v)$, $S$ is as in \cite{\QG, 23.1.6}; next, $S:\dUU_\ca@>>>\dUU_\ca$ is the retriction of $S:\dUU@>>>\dUU$
to $\dUU_\ca$; for general $A$, $S$ is obtained by applying $A\ot_\ca?$ to the $\ca$-linear 
map $S:\dUU_\ca@>>>\dUU_\ca$. 
Let $\o:\dUU_A@>>>\dUU_A$ be the algebra involution defined in \cite{\QG, 31.1.4}. 

\subhead 1.5\endsubhead
Let $\dBB$ be the canonical basis of $\dUU_\ca$ (see \cite{\QG, 25.2}). If $b\in\dBB$ we shall denote the element
$1\ot_\ca b\in A\ot_\ca\dUU_\ca=\dUU_A$ again by $b$. Thus $\dBB$ can be viewed also as an $A$-basis of $\dUU_A$.

Note that $1_\l\in\dBB$ for any $\l\in X$. We have $\dBB=\sqc_{\l,\l'\in X}(\dBB\cap(1_\l\dUU_A1_{\l'}))$. For 
$a,b$ in $\dBB$ we write $ab=\sum_{c\in\dBB}m_{a,b}^cc$ ($ab$ is the product in $\dUU_A$ and $m_{a,b}^c\in A$ is 
$0$ for all but finitely many $c$). For $a,b,c$ in $\dBB$ we define $\hm^{a,b}_c\in A$ by the following
requirement: for any $\l,\l',\l_1,\l'_1,\l_2,\l'_2$ in $X$ and any $c\in\dBB\cap(1_\l\dUU_A1_{\l'})$ we have 
$\D_{\l,\l',\l_1,\l'_1,\l_2,\l'_2}c=\sum_{a,b}\hm^{a,b}_ca\ot b$ where $a$ runs over 
$\dBB\cap(1_{\l_1}\dUU_A1_{\l'_1})$, $b$ runs over $\dBB\cap(1_{\l_2}\dUU_A1_{\l'_2})$ (in the last sum 
$\hm^{a,b}_c$ is $0$ for all but finitely many $(a,b)$).

For any $a\in\dBB$ we have $S(a)=s_a\ua$ where $a\m\ua$ is an involution of $\dBB$ and $s_a$ is $\pm$ a power of 
$v$ with $s_{\ua}=s_a$ (see \cite{\QG. 23.1.7, 26.3.2}). We have $s_{1_\l}=1$, $\un{1_\l}=1_{-\l}$ for $\l\in X$.
For any $a\in\dBB$ we have $\o(a)=e_aa^!$ where $a\m a^!$ is an involution of $\dBB$ and $e_a=\pm1$ (see 
\cite{\QG, 26.3.2}). We have $e_{1_\l}=1$, $(1_\l)^!=1_{-\l}$ for $\l\in X$.

Note that the quantities $m_{a,b}^c,\hm^{a,b}_c,s_a,e_a$ (in $A$) are the images of the corresponding quantities 
in $\ca$ under the given homomorphism $\ca@>>>A$ and $a\m\ua$, $a\m a^!$ are independent of $A$.

As in \cite{\QG, 25.4} for any $a,b,d,e$ in $\dBB$ we have

(i) $\sum_cm^c_{a,b}m^e_{c,d}=\sum_cm^e_{a,c}m^c_{b,d}$;

(ii) $\sum_c\hm^{a,b}_c\hm^{c,d}_e=\sum_c\hm^{a,c}_e\hm^{b,d}_c$;

(iii) $\sum_cm^c_{a,b}\hm^{e,d}_c=\sum_{a',b',c',d'}\hm^{a',b'}_a\hm^{c',d'}_bm^e_{a',c'}m^d_{b',d'}$;

(iv) $\hm^{a,b}_{1_\l}=1$ if $a=1_{\l'},b=1_{\l''}$, $\l'+\l''=\l$ and $\hm^{a,b}_{1_\l}=0$, otherwise.
\nl
From the definitions we have for any $a,b,c$ in $\dBB$:

(v) $s_cm_{a,b}^c=s_as_bm_{\ub,\ua}^{\uc}$;

(vi) $s_c\hm^{\ub,\ua}_{\uc}=s_as_b\hm^{a,b}_c$;

(vii) $e_cm_{a,b}^c=e_ae_bm_{a^!,b^!}^{c^!}$;

(viii) $e_c\hm^{b^!,a^!}_{c^!}=e_ae_b\hm^{a,b}_c$;

(ix) $\sum_{\l\in X}m_{a,1_\l}^c=\sum_{\l\in X}m_{1_\l,a}^c=\d_{a,c}$;

(x) $\hm^{a,1_0}_c=\hm^{1_0,a}_c=\d_{a,c}$;

(xi) $\sum_{d,e\in\dBB}\hm^{d,e}_am_{d,\te}^bs_e=\sum_{d,e\in\dBB}\hm^{d,e}_am_{\td,e}^bs_d$ is $1$ if
$b=1_\l,a=1_0$ and is $0$ otherwise;

(xii) $m_{a,b}^{1_0}=\d_{a,1_0}\d_{b,1_0}$;

(xiii) $\hm^{a,b}_c=\hm^{b,a}_c$ in $A$ if $v=1$ in $A$.

\subhead 1.6\endsubhead
Let $\fC_A$ be the category whose objects are $\dUU_A$-modules $M$ which are unital in the following sense: for 
any $z\in M$ we have $1_\l z=0$ for all but finitely many $\l\in X$ and $\sum_{\l\in X}1_\l z=z$ and which are 
finitely generated as an $A$-module. (A unital $\dUU_A$-module $M$ is also an $A$-module by 
$a:z\m \sum_{\l\in X}(a1_\l)z$ (all but finitely many terms of the last sum are $0$). A morphism in $\fC_A$ is a 
$\dUU_A$-linear map. When $A=\QQ(v)$ we write $\fC$ instead of $\fC_A$.

Let $M,M'\in\fC_A$, The tensor product $M\ot _AM'$ will be regarded as a $\dUU_A$-module by the rule 
$c(z\ot z')=\sum_{a,b\in\dBB}\hm^{a,b}_caz\ot bz'$ where $c\in\dBB$, $z\in M,z'\in M'$. (All but finitely many 
terms of the last sum are $0$.) The fact that the rule above defines an $\dUU_A$-module structure follows from 
1.5(iii). We have $M\ot_AM'\in\fC_A$. This makes $\fC_A$ into a monoidal tensor category.

For any $M,M'$ in $\fC_A$ let ${}_f\car_{M,M'}:M'\ot M@>>>M\ot M'$ be the isomorphism in $\fC_A$ defined in 
\cite{\QG, 32.1.5} in terms of a fixed function $f:X\T X@>>>\ZZ$ as in \cite{\QG, 32.1.3}. (In the formula for 
${}_f\car_{M,M'}$ in \cite{\QG, 32.1.4} we interpret $b^-m,b'{}^+m'$ as $\sum_\l b^-1_\l m$, 
$\sum_\l b'{}^+1_\l m'$.)

To any object $M\in\fC_A$ we associate a new object ${}^\o M\in\fC_A$ with the same underlying $A$-module as $M$ 
and such that for any $u\in\dUU_A$, the operator $u$ on ${}^\o M$ coincides with the operator $\o(u)$ on $M$. 

For any $\l\in X^+$ let $\L_\l=\bof/\ct_\l$ where $\ct_\l=\sum_i\bof\th_i^{\la i,\l\ra+1}$. Let $\et_\l=1+\ct_\l$. We regard $\L_\l$ as an object of $\fC$ in which 
$(\th_i^+1_{\l'})\et_\l=0$ for $i\in I,\l'\in X$, 
$(x^-1_{\l'})\et_\l=\d_{\l,\l'}x+\ct_\l$ for $x\in\bof,\l'\in X$.
We have $\L_\l\in\fC$ (\cite{\QG, 6.3.4}). 

Let $\L_{\l,\ca}=\dUU_\ca\et_\l$. This is an $\ca$-lattice in $\L_\l$ and 
$\L_{\l,\ca}\in\fC_\ca$. Note that $\et_\l\in\L_{\l,\ca}$.

Let $\L_{\l,A}=A\ot_\ca\L_{\l,\ca}$. By extension of scalars $\L_{\l,A}$ becomes a 
$\dUU_A$-module. We have
$\L_{\l,A}\in\fC_A$. We write $\et_\l$ instead of $1\ot_\ca\et_\l\in A\ot_\ca\L_{\l,\ca}=L_{\l,A}$. Now 
$\et_\l\in\L_{\l,A}$ regarded as an element of ${}^\o\L_{\l,A}$ will be denoted by $\x_{-\l}$. 

For $\l\in X^+$ there is a unique subset $\BB_\l$ of $\BB$ such that 
$b\m(b^-1_\l)\et_\l$ is a 
bijection of $\BB_\l$ onto a basis $\uBB_\l$ of $\L_{\l,A}$; for $b\in\BB-\BB_\l$ we 
have $(b^-1_\l)\et_\l=0$. There is a unique subset $\BB'_\l$ of $\BB$ such that 
$b\m(b^+1_{-\l})\x_{-\l}$ is a bijection of $\BB'_\l$ onto $\uBB_\l$; for 
$b\in\BB-\BB'_\l$ we have $(b^+1_{-\l})\x_{-\l}=0$. 

Note that $\BB_\l,\BB'_\l$ are independent of $A$.

\subhead 1.7\endsubhead
For $\l\in X^+$ let $\dUU^{\ge\l}$ be the set of all $u\in\dUU$ such that for any $\l'\in X^+$ such that 
$\l'\not\ge\l$ we have $u|_{\L_{\l'}}=0$ (a two-sided ideal of $\dUU$). By \cite{\QG, 29.1}, 
$\dUU^{\ge\l}\cap\dBB$ is a basis of $\dUU^{\ge\l}$ and there is a unique partition 
$\dBB=\sqc_{\l\in X^+}\dBB[\l]$ such that for any $\l\in X^+$ we have 
$\dUU^{\ge\l}\cap\dBB=\sqc_{\l'\in X^+;\l'\ge\l}\dBB[\l']$. By \cite{\QG, 29.1.6}, $\dBB[\l]$ is finite for any 
$\l\in X^+$.

\proclaim{Lemma 1.8} For any $c\in\dBB$, the set $\{(a,b)\in\dBB\T\dBB;m^c_{a,b}\in A-\{0\}\}$ is finite.
\endproclaim
We may assume that $A=\ca$. Let $c\in\dBB$. We have $c\in\dBB[\l]$ for a unique $\l\in X^+$. Let 
$(a,b)\in\dBB\T\dBB$ be such that $m^c_{a,b}\ne0$. We have $a\in\dBB[\l']$, $b\in\dBB[\l'']$ for some 
$\l',\l''$ in $X^+$. We have $a\in\dUU^{\ge\l'}$. Since 
$\dUU^{\ge\l'}\dUU\sub\dUU^{\ge\l'}$, we have 
$ab\in\dUU^{\ge\l'}$. Since $c$ appears with non-zero coefficient in $ab$ and $\dUU^{\ge\l'}\cap\dBB$ is a basis 
of $\dUU^{\ge\l'}$, we have $c\in\dUU^{\ge\l'}\cap\dBB$ that is $c\in\cup_{\l_1\in X^+;\l_1\ge\l'}\dBB[\l_1]$. Thus
$\l\ge\l'$. Similarly, $\l\ge\l''$. Thus $\{(a,b)\in\dBB\T\dBB;m^c_{a,b}\ne0\}$ is contained in
$\sqc_{(\l',\l'')\in X^+\T X^+;\l\ge\l',\l\ge\l''}\dBB[\l']\T\dBB[\l'']$; this is a finite set since 
$\{(\l',\l'')\in X^+\T X^+;\l\ge\l',\l\ge\l''\}$ is finite and each $\dBB[\l']\T\dBB[\l'']$ is finite. The lemma 
is proved.

\subhead 1.9\endsubhead
For any $\l,\l'$ in $X^+$ there is a unique (finite) subset $\dBB_{\l,\l'}$ of $\dBB$ such
that $a\m a(\x_{-\l}\ot\et_{\l'})$ defines a bijection of $\dBB_{\l,\l'}$ onto an
$A$-basis $\udBB_{\l,\l'}$ of 
${}^\o\L_{\l,A}\ot_A\L_{\l',A}$; moreover $a(\x_{-\l}\ot_A\et_{\l'})=0$ for any $a\in\dBB-\dBB_{\l,\l'}$. This
property appears in \cite{\QG, 25.2.1} in the case where $A=\QQ(v)$ but it is then automatically true for
general $A$. Note that the subset $\dBB_{\l,\l'}$ of $\dBB$ is independent of $A$. 
We have $\dBB=\cup_{\l,\l'\in X^+}\dBB_{\l,\l'}$. 

\subhead 1.10\endsubhead
We show:

(a) {\it If $M\in\fC_A$ then $\fA:=\{a\in\dBB;aM\ne0\}$ is finite.}
\nl
Let $z_1,z_2,\do,z_r$ be a set of generators of the $A$-module $M$ such that each $z_j$ is contained in $1_{\z_j}M$
for some $\z_j\in X$. As in the proof of \cite{\QG, 31.2.7}, for $j\in[1,r]$ we can find $\l_j,\l'_j\in X^+$ such
that $\l'_j-\l_j=\z_j$ and a morphism \lb
$\ph_j:{}^\o\L_{\l_j,A}\ot_A\L_{\l'_j,A}@>>>M$ in $\fC_A$ such that 
$\ph(\x_{-\l_j}\ot\et_{\l'_j})=z_j$. If $a\in\fA$ then for some $j$ we have $az_j\ne0$ 
hence 
$a(\x_{-\l_j}\ot\et_{\l'_j})\ne0$ in ${}^\o\L_{\l_j,A}\ot_A\L_{\l'_j,A}$.
Thus we have $a\in\dBB_{\l_j,\l'_j}$. We see that
$\fA\sub\cup_{j\in[1,r]}\dBB_{\l_j,\l'_j}$. Thus $\fA$ is finite.

\subhead 1.11\endsubhead
Let $\hUU_A$ be the $A$-module consisting of all formal linear combinations 
$$\sum_{a\in\dBB}n_aa$$ 
with $n_a\in A$. We define an $A$-algebra structure on $\hUU_A$ by 
$$(\sum_{a\in\dBB}n_aa)(\sum_{b\in\dBB}\tn_bb)=\sum_{c\in\dBB}r_cc$$
where $r_c=\sum_{(a,b)\in\dBB\T\dBB}m_{a,b}^cn_a\tn_b$. The last sum is well defined since by 1.8 it has only
finitely many non-zero terms. This algebra structure is associative by 1.5(i) and has a unit element 
$1=\sum_{\l\in X}1_\l$. We have an imbedding of $A$-algebras $\dUU_A\sub\hUU_A$ under which $b\in\dBB\sub\dUU_A$ 
coressponds to $\sum_{a\in\dBB}\d_{a,b}a\in\hUU_A$.

Let $\hUU^{(2)}_A$ be the $A$-module consisting of all formal linear combinations 
$$\sum_{(a,a')\in\dBB\T\dBB}n_{a,a'}a\ot a'$$
 with $n_{a,a'}\in A$. We define an $A$-algebra structure on $\hUU^{(2)}_A$ by
$$(\sum_{(a,a')\in\dBB\T\dBB}n_{a,a'}a\ot a')(\sum_{(b,b')\in\dBB\T\dBB}\tn_{b,b'}b\ot b')=
\sum_{(c,c')\in\dBB\T\dBB}r_{c,c'}c\ot c'$$
where $r_{c,c'}=\sum_{a,a',b,b'}m_{a,b}^cm_{a',b'}^{c'}n_{a,a'}\tn_{b,b'}$. (Again this sum is well defined by
1.8.) This algebra is associative by 1.5(i) with unit element $1=\sum_{\l,\l'\in X}1_\l\ot1_{\l'}$. Note that 
$\hUU^{(2)}_A$ is associated to the direct sum of two copies of our root datum in the same way as $\hUU_A$ is 
associated to our root datum. For $\x=\sum_{a\in\dBB}n_aa\in\hUU_A$, $\x'=\sum_{a'\in\dBB}\tn_{a'}a'\in\hUU_A$, we
set $\x\ot\x'=\sum_{(a,a')\in\dBB\T\dBB}n_a\tn_{a'}a\ot a'\in\hUU^{(2)}$. Note that $\x,\x'\m\x\ot\x'$ defines an
$A$-algebra homomorphisms $\hUU_A\ot_A\hUU_A@>>>\hUU_A^{(2)}$. We define an $A$-linear map 
$\D:\dUU_A@>>>\hUU^{(2)}_A$ by
$$c\m\sum_{(a,b)\in\dBB\T\dBB}\hm^{a,b}_ca\ot b$$
for any $c\in\dBB$. By 1.5(iii), $\D$ is an $A$-algebra homomorphism.

We define an $A$-linear map $\mm:\hUU^{(2)}_A@>>>\hUU_A$ by
$$\sum_{(a,a')\in\dBB\T\dBB}n_{a,a'}a\ot a'\m\sum_{b\in\dBB}(\sum_{(a,a')\in\dBB\T\dBB}n_{a,a'}m_{a,a'}^b)b.$$
This is well defined, by 1.8. In particular, if $\x,\x'\in\hUU_A$ we have $\mm(\x\ot\x')=\x\x'$ (product in 
$\hUU_A$) where $\x\ot\x'$ is regarded as an element of $\hUU^{(2)}$ as above.

We define an $A$-linear map $S:\hUU_A@>>>\hUU_A$ by $\sum_{a\in\dBB}n_aa\m\sum_{a\in\dBB}s_an_a\ua$ and $A$-linear maps $S^{(2)}:\hUU^{(2)}_A@>>>\hUU^{(2)}_A$,
$S^{(0,1)}:\hUU^{(2)}_A@>>>\hUU^{(2)}_A$ by  
$$S^{(2}:\sum_{(a,a')\in\dBB\T\dBB}n_{a,a'}a\ot a'\m\sum_{(a,a')\in\dBB\T\dBB}s_as_{a'}n_{a,a'}\ua\ot\ua',$$
$$S^{(0,1)}:\sum_{(a,a')\in\dBB\T\dBB}n_{a,a'}a\ot a'\m\sum_{(a,a')\in\dBB\T\dBB}s_{a'}n_{a,a'}a\ot\ua'.$$
In particular for $\x,\x'$ in $\hUU_A$ we have $S^{(2)}(\x\ot\x')=S(\x)\ot S(\x')$,
$S^{(0,1)}(\x\ot\x')=\x\ot S(\x')$.

We define an $A$-linear map $\o:\hUU_A@>>>\hUU_A$ by $\sum_{a\in\dBB}n_aa\m\sum_{a\in\dBB}e_an_aa^!$ and an 
$A$-linear map $\o^{(2)}:\hUU^{(2)}_A@>>>\hUU^{(2)}_A$ by 
$$\sum_{(a,a')\in\dBB\T\dBB}n_{a,a'}a\ot a'\m\sum_{(a,a')\in\dBB\T\dBB}e_ae_{a'}n_{a,a'}a^!\ot a'{}^!.$$
In particular for $\x,\x'$ in $\hUU_A$ we have $\o^{(2)}(\x\ot\x')=\o(\x)\ot \o(\x')$.

For $\x,\x'$ in $\hUU_A$ we have $S(\x\x')=S(\x')S(\x)$, $\o(\x\x')=\o(\x)\o(\x')$ (see 
1.5(v),(vii)).

We define an $A$-linear map $\e:\hUU_A@>>>A$ by $\sum_{a\in\dBB}n_aa\m n_{1_0}$. This is a homomorphism of
$A$-algebras preserving $1$, by 1.5(xii).

\subhead 1.12\endsubhead
Let $M\in\fC_A$. We define an $A$-linear map 
$\hUU_A@>>>\Hom_A(M,M)$ by $\sum_{a\in\dBB}n_aa:z\m\sum_{a\in\dBB}n_aaz$. The last sum
is well defined since, by 1.10(a), all but finitely many of its terms are $0$. This defines a structure of left
$\hUU_A$-module on $M$ extending the $\dUU_A$-module structure. Note that $1\in\hUU_A$ acts as $1$ on $M$.

Similarly, if $M\in\fC_A$, $M'\in\fC_A$ then the obvious $\dUU_A\ot\dUU_A$-module structure on $M\ot_AM'$ extends
to a $\hUU^{(2)}_A$-module structure on $M\ot_AM'$ defined by the $A$-linear map 
$\hUU^{(2)}_A@>>>\Hom_A(M\ot_AM',M\ot_AM')$,
$$\sum_{(a,a')\in\dBB\T\dBB}n_{a,a'}a\ot a':z\ot z'\m\sum_{(a,a')\in\dBB\T\dBB}n_{a,a'}az\ot a'z'.$$
 We can now 
regard $M\ot_AM'$ as an $\dUU_A$-module by restricting the $\hUU^{(2)}_A$-module structure to $\dUU_A$ via the 
algebra homomorphism $\D:\dUU_A@>>>\hUU^{(2)}_A$. Note that the resulting $\dUU_A$-module $M\ot_AM'$ is the one
defined in 1.6.

\subhead 1.13\endsubhead
Let $\l,\l',\l_1,\l_2,\l'_1,\l'_2$ be elements of $X^+$ such that $\l=\l_1+\l_2$, $\l'=\l'_1+\l'_2$. Let 
$\fT:\L_{\l'_1+\l'_2,\ca}@>>>\L_{\l'_1,\ca}\ot_\ca\L_{\l'_2,\ca}$ be the restriction of the unique unique morphism
$\L_{\l'_1+\l'_2}@>>>\L_{\l'_1}\ot\L_{\l'_2}$ in $\fC$ such that $\et_{\l'_1+\l'_2}\m\et_{\l'_1}\ot\et_{\l'_2}$, 
see \cite{\QG, 25.1.2(a), 25.1.2(b)}. Similarly let 
$\fT':{}^\o\L_{\l_1+\l_2,\ca}@>>>{}^\o\L_{\l_1,\ca}\ot_\ca{}^\o\L_{\l_2,\ca}$ be the restriction of the unique
morphism ${}^\o\L_{\l_1+\l_2}@>>>{}^\o\L_{\l_1}\ot{}^\o\L_{\l_2}$ in $\fC$ such that 
$\x_{-\l_1-\l_2}\m\x_{-\l_1}\ot\x_{-\l_2}$. Note that $\fT,\fT'$ are morphisms in $\fC_\ca$. We define 
$$\t:{}^\o\L_{\l,\ca}\ot_\ca\L_{\l',\ca}@>>>{}^\o\L_{\l_1,\ca}\ot_\ca\L_{\l'_1,\ca}\ot_\ca
{}^\o\L_{\l_2,\ca}\ot_\ca\L_{\l'_2,\ca}$$
as the composition
$$\align&{}^\o\L_{\l,\ca}\ot_\ca\L_{\l',\ca}@>\fT'\ot\fT>>
{}^\o\L_{\l_1,\ca}\ot_\ca{}^\o\L_{\l_2,\ca}\ot\L_{\l'_1,\ca}\ot_\ca\L_{\l'_2,\ca}\\&
@>1\ot v^{f(\l_2,\l'_1)}{}_f\car_{{}^\o\L_{\l_2,\ca},\L_{\l'_1,\ca}}\i\ot1>>
{}^\o\L_{\l_1,\ca}\ot_\ca\L_{\l'_1,\ca}\ot_\ca{}^\o\L_{\l_2,\ca}\ot_\ca\L_{\l'_2,\ca},\endalign$$
a morphism in $\fC_A$. 
Define $\r:\dUU_\ca@>>>{}^\o\L_{\l,\ca}\ot_\ca\L_{\l',\ca}$ by $u\m u(\x_{-\l}\ot\et_{\l'})$. Define 
$$\r':\hUU^{(2)}_\ca@>>>{}^\o\L_{\l_1,\ca}\ot_\ca\L_{\l'_1,\ca}\ot_\ca{}^\o\L_{\l_2,\ca}\ot_\ca\L_{\l'_2,\ca}$$
by $u\m u(\x_{-\l_1}\ot_\ca\et_{\l'_1}\ot_\ca\x_{-\l_2}\ot_\ca\et_{\l'_2})$. Here $M\ot_\ca M'$ with 
$M={}^\o\L_{\l_1,\ca}\ot_\ca\L_{\l'_1,\ca}$, $M'={}^\o\L_{\l_2,\ca}\ot_\ca\L_{\l'_2,\ca}$ is regarded as a
$\hUU^{(2)}_\ca$-module as in 1.12. We show that the diagram
$$\CD
\dUU_\ca@>\D>>\hUU^{(2)}_\ca\\
@V\r VV         @V\r'VV   \\
{}^\o\L_{\l,\ca}\ot_\ca\L_{\l',\ca}@>\t>>M\ot_\ca M'
\endCD$$
is commutative.
We regard $\dUU_\ca$ as a $\dUU_\ca$-module via left multiplication. We regard $\hUU^{(2)}_\ca$ as a 
$\dUU_\ca$-module in which $u\in\dUU_A$ acts by left multiplication by $\D(u)$. Then $\D$ is $\dUU_\ca$-linear 
since $\D$ is an algebra homomorphism. From the definitions we see that $\r$ and $\r'$ are $\dUU_\ca$-linear. Thus
all maps in our diagram are $\dUU_\ca$-linear. Since $\{1_x;x\in X\}$ generate the $\dUU_\ca$-module $\dUU_\ca$, it
is enough to show that the two compositions in the diagram coincide when applied to any $1_x$ with $x\in X$. Thus
it is enough to show that
$$\t(1_x(\x_{-\l}\ot\et_{\l'}))=   
\sum_{x',x''\in X;x'+x''=x}(1_{x'}(\x_{-\l_1}\ot\et_{\l'_1}))\ot(1_{x''}(\x_{-\l_2}\ot\et_{\l'_2}))$$
or equivalently that
$$\t(\x_{-\l}\ot\et_{\l'})=\x_{-\l_1}\ot\et_{\l'_1}\ot\x_{-\l_2}\ot\et_{\l'_2}.$$
This follows from the definitions using the equality
$${}_f\car_{{}^\o\L_{\l_2,\ca},\L_{\l'_1,\ca}}\i
(\x_{-\l_2}\ot\et_{\l'_1})=v^{-f(\l_2,\l'_1)}\et_{\l'_1}\ot\x_{-\l_2}.$$
We now show:

(a) {\it The $\ca$-linear map $\t$ is a split injection.}
\nl
Since $\t$ is the composition of $\fT'\ot_\ca\fT$ with an $\ca$-linear isomorphism, it is enough to show that 
$\fT'\ot_\ca\fT$ is a split injection. It is also enough to show that the $\ca$-linear maps $\fT$ and $\fT'$ are 
split injections. By \cite{\QG, 27.1.7}, $\fT$ carries $\uBB_{\l'}$ (an $\ca$-basis of $\L_{\l',\ca}$) bijectively
onto a subset of $\udBB_{\l'_1,\l'_2}$ (an $\ca$-basis of $\L_{\l'_1,\ca}\ot_\ca\L_{\l'_2,\ca}$). Hence $\fT$ is a split injection.

We identify the $\ca$-modules ${}^\o\L_{\l_1,\ca}\ot_\ca{}^\o\L_{\l_2,\ca}$, 
$\L_{\l_2,\ca}\ot_\ca\L_{\l_1,\ca}$ by $x\ot y\lra y\ot x$. Then $\fT'$ becomes the $\ca$-linear map 
$\L_{\l,\ca}@>>>\L_{\l_2,\ca}\ot_\ca\L_{\l_1,\ca}$ of the same type as $\fT$; hence it is a split injection by the
previous argument for $\fT$. This proves (a).

\subhead 1.14\endsubhead
Assume that $\l,\l',\l_1,\l_2,\l'_1,\l'_2$ in $X^+$ are such that $\l=\l_1+\l_2$, $\l'=\l'_1+\l'_2$. We show:

(a) {\it if $a\in\dBB_{\l_1,\l'_1}$, $b\in\dBB_{\l_2,\l'_2}$, $c'\in\dBB-\dBB_{\l,\l'}$ then $\hm^{a,b}_{c'}=0$ in
$\ca$.}
\nl
In the commutative diagram in 1.13 we have $\r(c')=0$ hence $\r'(\D(c'))=0$ that is
$\r'(\sum_{a,b\in\dBB}\hm^{a,b}_{c'}a\ot b)=0$ so that 
$$\sum_{a,b\in\dBB}\hm^{a,b}_{c'}(a(\x_{-\l_1}\ot\et_{\l'_1})\ot b(\x_{-\l_2}\ot\et_{\l'_2}))=0.$$
The term corresponding to $(a,b)$ is $0$ unless $(a,b)\in\dBB_{\l_1,\l'_1}\T\dBB_{\l_2,\l'_2}$; moreover when
$(a,b)$ runs through $\dBB_{\l_1,\l'_1}\T\dBB_{\l_2,\l'_2}$, the elements 
$a(\x_{-\l_1}\ot\et_{\l'_1})\ot b(\x_{-\l_2}\ot\et_{\l'_2})$ are linearly independent (they form the 
basis $\udBB_{\l_1,\l'_1}\ot\udBB_{\l_2,\l'_2}$). It follows that $\hm^{a,b}_{c'}=0$ for any 
$(a,b)\in\dBB_{\l_1,\l'_1}\T\dBB_{\l_2,\l'_2}$. This proves (a).

\subhead 1.15\endsubhead
In the setup of 1.14, we show:

(a) {\it Let $c\in\dBB_{\l,\l'}$. There exists a function $h:\dBB_{\l_1,\l'_1}\T\dBB_{\l_2,\l'_2}@>>>\ca$ such
that for any $c'\in\dBB$ we have}
$$\sum_{a\in\dBB_{\l_1,\l'_1},b\in\dBB_{\l_2,\l'_2}}h(a,b)\hm^{a,b}_{c'}=\d_{c,c'}.$$
For any $\ti\l,\ti\l'$ in $X^+$, we set $b^\da=b(\x_{-\ti\l}\ot_\ca\et_{\ti\l'})$ for $b\in\dBB_{\ti\l,\ti\l'}$ so
that $\{b^\da;b\in\dBB_{\ti\l,\ti\l'}\}$ is an $\ca$-basis of ${}^\o\L_{\ti\l,\ca}\ot_\ca\L_{\ti\l',\ca}$ and we
denote by $\{\hb;b\in\dBB_{\ti\l,\ti\l'}\}$ the dual $\ca$-basis of 
$\Hom_\ca({}^\o\L_{\ti\l,\ca}\ot_\ca\L_{\ti\l',\ca},\ca)$. In particular the bases $\{\hc';c'\in\dBB_{\l,\l'}\}$,
$\{\ha';a'\in\dBB_{\l_1,\l'_1}\}$, $\{\hb';b'\in\dBB_{\l_2,\l'_2}\}$ of 
$\Hom_\ca(\L_{\l,\ca}\ot_\ca\L_{\l',\ca},\ca)$, $\Hom_\ca(\L_{\l_1,\ca}\ot_\ca\L_{\l'_1,\ca},\ca)$,
$\Hom_\ca(\L_{\l_2,\ca}\ot_\ca\L_{\l'_2,\ca},\ca)$, are defined.
By 1.13(a), we can find an $\ca$-linear map
$$\ps:\L_{\l_1,\ca}\ot_\ca\L_{\l'_1,\ca}\ot_\ca\L_{\l_2,\ca}\ot_\ca\L_{\l'_2,\ca}@>>>
\L_{\l,\ca}\ot_\ca\L_{\l',\ca}$$
such that $\ps\t=1$ with $\t$ as in 1.13. The transpose maps ${}^t\ps,{}^t\t$ satisfy ${}^t\t({}^t\ps)=1$. We can 
write uniquely 
$${}^t\ps(\hc)=\sum_{a'\in\dBB_{\l_1,\l'_1},b'\in\dBB_{\l_2,\l'_2}}h(a',b')\ha'\ot\hb'$$
with $h(a',b')\in\ca$. Hence 
$$\hc=\sum_{a'\in\dBB_{\l_1,\l'_1},b'\in\dBB_{\l_2,\l'_2}}h(a',b'){}^t\t(\ha'\ot\hb').$$
Evaluating both sides at $c'{}^\da$ with $c'\in\dBB_{\l,\l'}$ and using
$$\t(c'{}^\da)=\sum_{a\in\dBB_{\l_1,\l'_1},b\in\dBB_{\l_2,\l'_2}}\hm^{a,b}_{c'}a^\da\ot b^\da$$
(which follows from the commutative diagram in 1.13) we have
$$\d_{c,c'}=\sum_{a,a'\in\dBB_{\l_1,\l'_1},b,b'\in\dBB_{\l_2,\l'_2}}h(a',b')\hm^{a,b}_{c'}\d_{a',a}\d_{b',b}$$
that is
$$\d_{c,c'}=\sum_{a\in\dBB_{\l_1,\l'_1},b\in\dBB_{\l_2,\l'_2}}h(a,b)\hm^{a,b}_{c'}.$$
This also holds when $c'\in\dBB-\dBB_{\l,\l'}$ (both sides are $0$ by 1.14(a)). This proves (a).

\proclaim{Lemma 1.16} For any $a,b$ in $\dBB$, the set $\{c\in\dBB;\hm^{a,b}_c\in\ca-\{0\}\}$ is finite.
\endproclaim
We have $a\in\dBB_{\l_1,\l'_1}$, $b\in\dBB_{\l_2,\l'_2}$ for some 
$\l_1,\l'_1,\l_2,\l'_2$ in $X^+$. By 1.14(a) we have 
$\{c\in\dBB;\hm^{a,b}_c\in\ca-\{0\}\}\sub\dBB_{\l_1+\l_2,\l'_1+\l'_2}$. Since $\dBB_{\l_1+\l_2,\l'_1+\l'_2}$ is 
finite, the lemma holds.

\subhead 1.17\endsubhead
We define an $A$-linear map $\hD:\hUU_A@>>>\hUU^{(2)}_A$ by
$$\sum_{c\in\dBB}n_cc\m\sum_{(a,b)\in\dBB\T\dBB}(\sum_{c\in\dBB}\hm^{a,b}_cn_c)a\ot b.$$
This make sense: for any $(a,b)\in\dBB\T\dBB$, the sum $\sum_{c\in\dBB}\hm^{a,b}_cn_c$ has only finitely many 
non-zero terms. (See 1.16.) Using 1.5(iii) we see that $\hD$ is an $A$-algebra homomorphism. It clearly extends 
the homomorphism $\D:\dUU_A@>>>\hUU^{(2)}_A$. From 1.5(vi),(viii) we see that 
$$\hD(S(\x))=S^{(2)}\text{\rm tr}\hD(\x),\qua \hD(\o(\x))=\o^{(2)}\text{\rm tr}\hD(\x)$$
for any $\x\in\hUU_A$ where $S^{(2)},\o^{(2)}$ are as in 1.11 and $\text{\rm tr}:\hUU^{(2)}_A@>>>\hUU^{(2)}_A$ is
the $A$-linear map given by 
$$\sum_{(a,a')\in\dBB\T\dBB}n_{a,a'}a\ot a'\m\sum_{(a,a')\in\dBB\T\dBB}n_{a',a}a\ot a'.$$
If $v=1$ in $A$ then $\text{\rm tr}\hD(\x)=\hD(\x)$ for $\x\in\hUU_A$ (see 1.5(xiii)). Hence in this case
we have $\hD(S(\x))=S^{(2)}\hD(\x)$, $\hD(\o(\x))=\o^{(2)}\hD(\x)$ for $\x\in\hUU_A$.

\subhead 1.18\endsubhead
For any $i\in I$ and $h\in A$ we set
$$x_i(h)=\sum_{c\in\NN,\l\in X}h^c\th_i^{(c)+}1_\l\in\hUU_A,$$
$$y_i(h)=\sum_{c\in\NN,\l\in X}h^c\th_i^{(c)-}1_\l\in\hUU_A.$$
Note that $\th_i^{(t)\pm}1_\l\in\dBB$, see \cite{\QG, 25.3.1}. We have
$$y_i(h)=\o(x_i(h)).$$
We show:

(a) {\it If $v=1$ in $A$ then} $\hD(x_i(h))=x_i(h)\ot x_i(h)\in\hUU^{(2)}_A$, 
$\hD(y_i(h))=y_i(h)\ot y_i(h)\in\hUU^{(2)}_A$.
\nl
Without assumption on $v\in A$ we have
$$\hD(x_i(h))=\sum_{t',t''\in\NN,\l',\l''\in X}h^{t'+t''}v_i^{t't''}v_i^{<i,\l'>t''}(\th_i^{(t')+}1_{\l'})
\ot(\th_i^{(t''+)}1_{\l''}),$$
$$\hD(y_i(h))=\sum_{t',t''\in\NN,\l',\l''\in X}h^{t'+t''}v_i^{-t't''}v_i^{-<i,\l''>t'}(\th_i^{(t')-}1_{\l'})
\ot(\th_i^{(t'')-}1_{\l''}).$$
The desired formulas follow.

We show:

(b) {\it If $v=1$ in $A$ and $h,h'\in A$, then} $x_i(h+h')=x_i(h)x_i(h')\in\hUU_A$, 
$y_i(h+h')=y_i(h)y_i(h')\in\hUU_A$.
\nl
Without assumption on $v\in A$ we have
$$x_i(h)x_i(h')=\sum_{t\in\NN,\l\in X}(\sum_{t',t'';t'+t''=t}h^{t'}h'{}^{t''}\bcc{t}{t'}_i)\th_i^{(t)+}1_\l,$$
$$y_i(h)y_i(h')=\sum_{t\in\NN,\l\in X}(\sum_{t',t'';t'+t''=t}h^{t'}h'{}^{t''}\bcc{t}{t'}_i)\th_i^{(t)-}1_\l.$$
If $v=1$ we have $\bcc{t}{t'}_i=\bin{t}{t'}$ and $\sum_{t',t'';t'+t''=t}h^{t'}h'{}^{t''}\bcc{t}{t'}_i=(h+h')^t$.
The desired formulas follow.

\subhead 1.19\endsubhead
We say that a (possibly infinite) sum $\sum_{f\in F}p_f$ with $p_f\in\hUU_A$ is {\it defined in $\hUU_A$} if, 
when setting $p_f=\sum_{a\in\dBB}n_{a,f}a\in\hUU_A$ ($n_{a,f}\in A$), the set $\{f\in F;n_{a,f}\ne0\}$ is finite 
for any $a\in\dBB$. In this case we assign to the sum the value $\sum_{a\in\dBB}(\sum_{f\in F}n_{a,f})a\in\hUU_A$.
Using 1.8 we see that, if the (possibly infinite) sums $\sum_{f\in F}p_f$, $\sum_{f'\in F'}p'_{f'}$ are defined in
$\hUU_A$, then the sum $\sum_{(f,f')\in F\T F'}p_fp'_{f'}$ is defined in $\hUU_A$ and we have
$$(\sum_{f\in F}p_f)(\sum_{f'\in F'}p'_{f'})=\sum_{(f,f')\in F\T F'}p_fp'_{f'},\text { in }
\hUU_A.$$

\head 2. The elements $s'_{i,e},s''_{i,e}$ of $\hUU_A$\endhead
\subhead 2.1\endsubhead
In this section we introduce and study some elements
$s'_{i,e},s''_{i,e}$ of $\hUU_A$ which implement the braid group action \cite{\QG, Ch.39}
on $\dUU_A$.

Let $i\in I,e=\pm1$. Let
$$s'_{i,e}=\sum_{y\in\ZZ,x\in\ZZ,\l\in X;\la i,\l\ra=x+y}(-1)^yv_i^{ey}\th_i^{(x)-}1_\l\th_i^{(y)+}\in\hUU_A,$$
$$s''_{i,e}=\sum_{y\in\ZZ,x\in\ZZ,\l\in X;\la i,\l\ra=x+y}(-1)^xv_i^{ex}\th_i^{(x)-}1_\l\th_i^{(y)+}\in\hUU_A.$$
Note that in the formulas above each $\th_i^{(x)-}1_\l\th_i^{(y)+}$ belongs to $\dBB$ (see \cite{\QG, 25.3.1}). We
have $\o(s'_{i,e})=s''_{i,e}$.

\proclaim{Lemma 2.2} For any $r\in\ZZ$ we consider the sums
$$\sum_{m,n,p\in\NN,\l\in X;\la i,\l\ra=m-n+p-r}(-1)^nv_i^{e(-mp+n)}\th_i^{(m)-}\th_i^{(n)+}\th_i^{(p)-}1_\l,
\tag i$$
$$\sum_{m,n,p\in\NN,\l\in X;\la i,\l\ra=-m+n-p+r}(-1)^nv_i^{e(-mp+n)}\th_i^{(m)+}\th_i^{(n)-}\th_i^{(p)+}1_\l.
\tag ii$$
(a) The sums (i), (ii) are defined in $\hUU_A$. Let $\t'_{i,e,r},\t''_{i,e,r}$ be their value in $\hUU_A$.

(b) We have $\t'_{i,e,0}=s'_{i,e}$, $\t''_{i,e,0}=s''_{i,e}$.

(c) If $v=1$ in $A$ and $r\in\ZZ-\{0\}$, we have $\t'_{i,e,r}=0$, $\t''_{i,e,r}=0$.

(d) If $v=1$ in $A$ we have $s''_{i,e}=x_i(1)y_i(-1)x_i(1)$, $s'_{i,e}=y_i(1)x_i(-1)y_i(1)$ in $\hUU_A$.

(e) If $v=1$ in $A$ we have $\hD(s'_{i,e})=s'_{i,e}\ot s'_{i,e}$ and $\hD(s''_{i,e})=s''_{i,e}\ot s''_{i,e}$ in
$\hUU^{(2)}_A$.
\endproclaim
We compute formally the sum (ii) for $r\in\NN$ (using 1.4):
$$\align&\sum_{m,n,p\in\NN,\l\in X;\la i,\l\ra=-m+n-p+r}
(-1)^nv_i^{e(-mp+n)}\th_i^{(m)+}1_{\l+pi'-ni'}\th_i^{(n)-}\th_i^{(p)+}
\\&=\sum_{m\in\ZZ,n\in\ZZ,p\in\NN,t\in\NN,\l\in X;\la i,\l\ra=-m+n-p+r}
(-1)^nv_i^{e(-mp+n)}\T\\&\T\bcc{m+n-\la i,-\l-pi'+ni'\ra}{t}_i
\th_i^{(n-t)-}1_{\l+pi'-ni'+(n+m-t)i'}\th_i^{(m-t)+}\th_i^{(p)+}
\\&=\sum_{m\in\ZZ,n\in\ZZ,p\in\NN,t\in\NN,\l\in X;\la i,\l\ra=-m+n-p+r}
(-1)^nv_i^{e(-mp+n)}\bcc{r+p}{t}_i\T\\&\T\bcc{m+p-t}{p}_i\th_i^{(n-t)-}1_{\l+pi'-ni'+(n+m-t)i'}\th_i^{(m+p-t)+}
\\&=\sum_{m\in\ZZ,x\in\ZZ,p\in\NN,t\in\NN,\l\in X;\la i,\l\ra=-m+x+t-p+r}
(-1)^{x+t}v_i^{e(-mp+x+t)}\bcc{r+p}{t}_i\T\\&\T
\bcc{m+p-t}{p}_i\th_i^{(x)-}1_{\l+(m+p-t)i'}\th_i^{(m+p-t)+}
\\&=\sum\Sb y\in\NN,x\in\NN,\l\in X;\\ \la i,\l\ra=-y+x+r\eSb
(-1)^{x+t}v_i^{e(-(y-p+t)p+x+t)}\bcc{r+p}{t}_i\bcc{y}{p}_i
c_{x,y,\l;r}\th_i^{(x)-}1_{\l+yi'}\th_i^{(y)+}\endalign$$
where the last product $\th_i^{(x)-}1_{\l+yi'}\th_i^{(y)+}$ belongs to $\dBB$ (see \cite{\QG, 25.3.1}) and
$$c_{x,y,\l;r}=\sum_{p\in\NN,t\in\NN}(-1)^{x+t}v_i^{e(-(y-p+t)p+x+t)}\bcc{r+p}{t}_i\bcc{y}{p}_i.$$
Note that only finitely many terms of this sum can be $\ne0$ (those for which $p\le y$ and $t\le r+p$); here we 
have used that $r\in\NN$. In particular the sum (ii) is defined in $\hUU_A$ and has a value $\t''_{i,e;r}$ in
$\hUU_A$.

Next we compute formally the sum (ii) for $r\in-\NN$ (using 1.4):
$$\align&\sum_{m,n,p\in\NN,\l\in X;\la i,\l\ra=-m+n-p+r}
(-1)^nv_i^{e(-mp+n)}\th_i^{(m)+}\th_i^{(n)-}1_{\l+pi'}\th_i^{(p)+}\\&
=\sum_{m\in\NN,n\in\ZZ,p\in\ZZ,t\in\NN,\l\in X;\la i,\l\ra=-m+n-p+r}
(-1)^nv_i^{e(-mp+n)}\bcc{n+p-\la i,\l+pi'\ra}{t}_i\T\\&\T\th_i^{(m)+}\th_i^{(p-t)+}1_{\l+pi'-(n+p-t)i'}\th_i^{(n-t)-}
\\&=\sum_{m\in\NN,n\in\ZZ,p\in\ZZ,t\in\NN,\l\in X;\la i,\l\ra=-m+n-p+r}
(-1)^nv_i^{e(-mp+n)}\bcc{m-r}{t}_i\T\\&\T\bcc{m+p-t}{m}_i\th_i^{(m+p-t)+}1_{\l-(n-t)i'}\th_i^{(n-t)-}\\&
=\sum_{m\in\NN,x\in\ZZ,y\in\ZZ,t\in\NN,\l\in X;\la i,\l\ra=y-x+r}
(-1)^{x+t}v_i^{e(-m(x+t-m)+y+t)}\bcc{m-r}{t}_i\T\\&\T\bcc{x}{m}_i\th_i^{(x)+}1_{\l-yi'}\th_i^{(y)-}\\&
=\sum_{x\in\NN,y\in\NN,\l\in X;\la i,\l\ra=y-x+r}\tc_{x,y,\l;r}\th_i^{(x)+}1_{\l-yi'}\th_i^{(y)-},\endalign$$
where the last product $\th_i^{(x)+}1_{\l-yi'}\th_i^{(y)-}$ belongs to $\dBB$ (see \cite{\QG, 25.3.1}) and
$$\tc_{x,y,\l;r}=\sum_{m\in\NN,t\in\NN}(-1)^{x+t}v_i^{e(-m(x+t-m)+y+t)}\bcc{m-r}{t}_i\bcc{x}{m}_i.$$
Note that only finitely many terms of this sum can be $\ne0$ (those for which $m\le x$ and $t\le m-r$); here we 
have used that $r\in-\NN$. In particular the sum (ii) is defined in $\hUU_A$ and has a value $\t''_{i,e;r}$ in
$\hUU_A$.

From the computations above we see also that the sum 
$$\sum_{m,n,p\in\NN,\l\in X}(-1)^nv_i^{e(-mp+n)}\th_i^{(m)+}\th_i^{(n)-}\th_i^{(p)+}1_\l$$
is defined in $\hUU_A$ and its value $\t''_{i,e}\in\hUU_A$ is equal to $\sum_{r\in\ZZ}\t''_{i,e,r}$.

In the formula for $c_{x,y,\l;0}$ we have
$$\sum_{t\in\NN}(-1)^{x+t}v_i^{e(-(y-p+t)p+x+t)}\bcc{r+p}{t}_i=
(-1)^xv_i^{e(-(y-p)p+x)}\d_{p,0}=(-1)^xv_i^{ex}\d_{p,0}$$
hence $c_{x,y,\l;0}=(-1)^xv_i^{ex}$ and  
$$\t''_{i,e;0}=\sum_{y\in\ZZ,x\in\ZZ,\l\in X;\la i,\l\ra=-y+x}(-1)^xv_i^{ex}\th_i^{(x)-}1_{\l+yi'}\th_i^{(y)+}.$$
We see that $\t''_{i,e;0}=s''_{i,e}$. Similarly, we see that the family defining the sum (i) is defined hence it 
has a value $\t'_{i,e;r}$ in $\hUU_A$, that the sum 
$$\sum_{m,n,p\in\NN,\l\in X}(-1)^nv_i^{e(-mp+n)}\th_i^{(m)-}\th_i^{(n)+}\th_i^{(p)-}1_\l,$$
is defined in $\hUU_A$ and its value $\t'_{i,e}\in\hUU_A$ is equal to $\sum_{r\in\ZZ}\t'_{i,e,r}$. We see also
that $\t'_{i,e;0}=s'_{i,e}$. (These statements could be also deduced from the earlier part of the proof using
that $\t''_{i,e;r}=\o(\t'_{i,e.r})$.)

In the remainder of the proof we assume that $v=1$ in $A$. In this case in the formula defining $c_{x,y,\l;r}$ 
with $r>0$ we have $\sum_{t\in\NN}(-1)^{x+t}\bcc{r+p}{t}_i=0$ since $r+p>0$ and $v=1$. Hence $c_{x,y,\l;r}=0$ and
$\t''_{i,e;r}=0$ for $r>0$. In the formula defining $\tc_{x,y,\l;r}$ with $r<0$ we have 
$\sum_{t\in\NN}(-1)^{x+t}\bcc{m-r}{t}_i=0$ since $m-r>0$ hence $\tc_{x,y,\l;r}=0$ and $\t''_{i,e;r}=0$ for $r<0$.
We see that $\t''_{i,e;r}=0$ for $r\ne0$ and $\t''_{i,e}=\t''_{i,e;0}=s''_{i,e}$. Applying $\o$ we see that 
$\t'_{i,e;r}=0$ for $r\ne0$ and $\t'_{i,e}=\t'_{i,e;0}=s'_{i,e}$.
It is clear that $\t''_{i,e}=x_i(1)y_i(-1)x_i(1)$, $\t'_{i,e}=y_i(1)x_i(-1)y_i(1)$ in $\hUU_A$. Hence (d) follows.
Now (e) follows from (d) using 1.18(a) and the fact that $\hD:\hUU_A@>>>\hUU^{(2)}_A$ is an algebra homomorphism.

\subhead 2.3\endsubhead
Let $i\in I,e=\pm1$. Let $T'_{i,e},T''_{i,e}$ be the algebra automorphisms of $\dUU_A$ defined in
\cite{\QG, 41.1.8}. For any $M\in\fC_A$ let $T'_{i,e}:M@>>>M$, $T''_{i,e}:M@>>>M$ be the $A$-linear isomorphisms
defined in \cite{\QG, 41.2.3}. From the definitions we have $T'_{i,e}(z)=\t'_{i,e;0}z$, 
$T''_{i,e}(z)=\t''_{i,e;0}z$ for $z\in M$ (using the $\hUU_A$-module structure of $M$). By \cite{\QG, 41.2.4}, for
$u\in\dUU_A,z\in M$ we have 
$$T'_{i,e}(uz)=T'_{i,e}(u)(T'_{i,e}(z)),\qua T''_{i,e}(uz)=T''_{i,e}(u)(T''_{i,e}(z)),$$
hence, using 2.2(b):
$$s'_{i,e}uz=T'_{i,e}(u)(s'_{i,e}z),\qua s''_{i,e}uz=T''_{i,e}(u)(s''_{i,e}z).\tag a$$
We show:
$$s'_{i,e}s''_{i,-e}=1.\tag b$$
Let $s=s'_{i,e}s''_{i,-e}-1\in\hUU_A$. We have $s=\sum_{a\in\dBB}n_aa$ with $n_a\in A$. If $M\in\fC_A$ and 
$z\in M$ we have $sz=T'_{i,e}T''_{i,-e}(z)-z=0$, see \cite{\QG, 41.2.4}. In particular if $\l,\l'\in X^+$ we have 
$s(\xi_{-\l}\ot\et_\l)=0$ in ${}^\o\L_{\l,A}\ot_A\L_{\l',A}$. Hence 
$\sum_{a\in\dBB_{\l,\l'}}n_aa(\xi_{-\l}\ot\et_\l)=0$. Since the elements $a(\xi_{-\l}\ot\et_\l)$ with 
$a\in\dBB_{\l,\l'}$ form a basis of ${}^\o\L_{\l,A}\ot_A\L_{\l',A}$ it follows that $n_a=0$ for any 
$a\in\dBB_{\l,\l'}$. Since $\dBB=\cup_{\l,\l'}\dBB_{\l,\l'}$ we see that $n_a=0$ for any $a\in\dBB$ hence $s=0$
and (b) is proved.

We show:

(c) {\it for any $u\in\dUU_A$ we have $T'_{i,e}(u)=s'_{i,e}us'_{i,e}{}\i$, $T''_{i,e}(u)=s''_{i,e}us''_{i,e}{}\i$.
}\nl
If $M\in\fC_A$ and $z\in M$ we have (using (b),(a)) $s'_{i,e}us'_{i,e}{}\i s'_{i,e}z=s'_{i,e}uz=T'_{i,e}(u)(z)$. 
Thus, setting $u'=s'_{i,e}us'_{i,e}{}\i-T'_{i,e}(u)\in\dUU_A$, we have $u'z=0$. Since this holds for any $M,z$ we 
see as in the proof of (b) that $u'=0$. This proves the first equality in (c). The second equality is proved
similarly.

\subhead 2.4\endsubhead
Let $i\ne j$ in $I$. Let $e=\pm1$. Let $n=n_{ij}=n_{ji}$ be as in 1.1. We show:

(a) $s'_{i,e}s'_{j,e}s'_{i,e}\do=s'_{j,e}s'_{i,e}s'_{j,e}\do$ in $\hUU_A$;

(b) $s''_{i,e}s''_{j,e}s''_{i,e}\do=s''_{j,e}s''_{i,e}s''_{j,e}\do$ in $\hUU_A$;
\nl
(all products in (a) and (b) have $n$ factors).
\nl
Let $s_1\in\hUU_A$ be the left hand side minus the right hand side of (a). If $M\in\fC_A$ and $z\in M$, we have 
$s_1z=(T'_{i,e}T'_{j,e}T'_{i,e}\do)z-(T'_{j,e}T'_{i,e}T'_{j,e}\do)z$ and this is $0$ by \cite{\QG, 41.2.4(a)}.
Thus $s_1z=0$. Since this holds for any $M,z$ we see as in the proof of 2.3(b) that $s_1=0$. Thus (a) holds. The
proof of (b) is similar.

\subhead 2.5\endsubhead
Let $e=\pm1$. Let $w\in W$. From 2.4(a),(b) we deduce by a standard argument that there are unique elements 
$w'_e\in\hUU_A$, $w''_e\in\hUU_A$ such that $w'_e=s'_{i_1,e}s'_{i_2,e}\do s'_{i_r,e}$, 
$w''_e=s''_{i_1,e}s''_{i_2,e}\do s''_{i_r,e}$ for any sequence $i_1,i_2,\do i_r$ in $W$ such that
$w=s_{i_1}s_{i_2}\do s_{i_r}$ in $W$ with $r=l(w)$. Using 2.3(b) we see that $w'_e(w\i)''_{-e}=1$. From the 
definitions we have $\o(w'_e)=w''_e$.

\subhead 2.6\endsubhead
Let $e=\pm1$. Assume that $v=1$ in $A$. From 2.2(d) we see that for $i\in I$, $s'_{i,e},s''_{i,e}$ are independent
of the choice of $e$; we denote them by $s'_i,s''_i$. Using this and the definitions we see that for $w\in W$, 
$w'_e,w''_e$ are independent of the choice of $e$; we denote them by $w',w''$. Using 2.2(e) and the fact 
that $\hD$ is an algebra homomorphism we see that

(a) $\hD(w')=w'\ot w'$, $\hD(w'')=w''\ot w''$ (in $\hUU^{(2)}_A$) for any $w\in W$.

\subhead 2.7\endsubhead
Let $\dBB^+=\{b^+1_\l;b\in\BB,\l\in X\}=\{1_{\l'}b^+;b\in\BB,\l'\in X\}$, a subset of $\dBB$, see 
\cite{\QG, 25.2.6}. Let $\dBB^-=\{b^-1_\l;b\in\BB,\l\in X\}=\{1_{\l'}b^-;b\in\BB,\l'\in X\}$, a subset of $\dBB$.
Let $\dUU^+_A$ (resp. $\dUU^-_A$) be the $A$-submodule of $\dUU_A$ spanned by $\dBB^+$ (resp. by $\dBB^-$). Now 
$\dBB^+$ (resp. $\dBB^-$) is a basis of the $A$-module $\dUU^+_A$ (resp. $\dUU^-_A$) and $\dUU^+_A,\dUU^-_A$ are 
subalgebras of $\dUU_A$. From the definitions we have $\o(b^+1_\l)=b^-1_{-\l}$ for $b\in\BB$, $\l\in X$. Hence
setting $a=b^+1_\l\in\dBB^+$ we have $e_a=1$ and $a^!=b^-1_{-\l}$ (see 1.5). Thus, $\o:\dUU_A@>>>\dUU_A$ restricts
to a bijection $\dBB^+@>\si>>\dBB^-$ and to an $A$-algebra isomorphism $\dUU^+_A@>\si>>\dUU^-_A$. 

For $b,b'\in\BB$ we set $b^+=\sum_{\l\in X}b^+1_\l\in\hUU_A$, $b'{}^-=\sum_{\l\in X}b'{}^-1_\l\in\hUU_A$. Then the
elements $b^+1_\l b'{}^-\in\dUU_A$, $b'{}^-1_\l b^-\in\dUU_A$ (as in 1.4) which are not products in $\dUU_A$, can
be interpreted as products $b^+\cdot 1_\l\cdot b'{}^-$, $b'{}^-\cdot 1_\l\cdot b^+$ in $\hUU_A$.

\subhead 2.8\endsubhead
Assume that $v=1$ in $A$. Let $w_0$ be as in 1.1. 
Let $n=l(w_0)$. We fix a sequence $i_1,i_2,\do,i_n$ in $I$ such that $s_{i_1}s_{i_2}\do s_{i_n}=w_0$. For any 
$c\in\NN$ and $k\in[1,n]$ there is a unique element $x_{c,k}\in\bof_A$ such that in $\dUU_A$ we have
$$T''_{i_1,1}T''_{i_2,1}\do T''_{i_{k-1},1}(\th_{i_k}^{(c)+}1_{s_{i_{k-1}}s_{i_{k-2}}\do s_{i_1}\l})=x_{c,k}^+1_\l
$$
for any $\l\in X^+$. (See \cite{\QG, 41.1.3}.) For any $\boc=(c_1,c_2,\do,c_n)\in\NN^n$ we set
$$x_\boc=x_{c_1,1}x_{c_2,2}\do x_{c_n,n}\in\bof_A.$$
By \cite{\QG, 41.1.4}, \cite{\QG, 41.1.7}, the set $\{x_\boc;\boc\in\NN^n\}$ is an $A$-basis of $\bof_A$. We show:

(a) {\it For any $\boc\in\NN^n$ we have
$$\hD(x_\boc^+)=
\sum_{\boc',\boc''\in\NN;\boc'+\boc''=\boc}x_{\boc'}^+\ot x_{\boc''}^+\text{ in }\hUU^{(2)}.$$}
(The last sum is defined since $\boc',\boc''$ only take finitely many values.)
For any $k\in[1,n],c\in\NN$ we have
$$\hD(x_{c,k}^+)=\sum\Sb c',c''\in\NN\\c'+c''=c\eSb x_{c',k}^+\ot x_{c'',k}^+.$$
(We use the formulas $T''_{i,1}(u)=s''_{i,1}us''_{i,1}{}\i$, see 2.3(c), and 
$\hD(s''_{i,1})=s''_{i,1}\ot s''_{i,1}$, see 2.2(e).) It follows that 
$$\align&\hD(x_\boc^+)=
\hD(x_{c_1,1}^+)\hD(x_{c_2,2}^+)\do\hD(x_{c_n,n}^+)\\&
=\sum_{\boc',\boc''\in\NN;\boc'+\boc''=\boc}
(x_{c'_1,1}^+x_{c'_2,2}^+\do x_{c'_n,n}^+)\ot
(x_{c''_1,1}^+x_{c''_2,2}^+\do x_{c''_n,n}^+)\endalign$$
which yields (a).

\head 3. The Hopf algebra $\OO_A$\endhead
\subhead 3.1\endsubhead
In this section we define the Hopf algebra $\OO_A$ as a submodule of the dual of 
$\dUU_A$ defined in terms of $\dBB$. We also study some basis properties of $\OO_A$.

For any $A$-module $V$ we set $V^\di=\Hom_A(V,A)$. For any $a\in\dBB$ we define a linear form $a^*:\dUU_A@>>>A$ by
$a'\m\d_{a,a'}$ for all $a'\in\dBB$. Let $\OO_A$ be the $A$-submodule of $\dUU_A^\di$ spanned by 
$\{a^*;a\in\dBB\}$. Thus $\{a^*;a\in\dBB\}$ is an $A$-basis of $\OO_A$. We define an $A$-algebra structure on 
$\OO_A$ by the rule $a^*b^*=\sum_{c\in\dBB}\hm^{a,b}_cc^*$ for any $a,b$ in $\dBB$. (The sum is well defined by 
1.16.) This algebra structure has a unit element namely $1_0^*$. The $A$-linear map $\d:\OO_A@>>>\OO_A\ot_A\OO_A$
given by $c^*\m\sum_{(a,b)\in\dBB\T\dBB}m_{a,b}^ca^*\ot b^*$ is well defined by 1.8; we call it comultiplication.
Define an $A$-linear map (antipode) $S:\OO_A@>>>\OO_A$ by $S(a^*)=s_a\ua^*$ (see 1.5) for $a\in\dBB$. 
Using 1.5(i)-(vi), 1.5(ix)-(xii) we see that $(\OO_A,\d,S)$ is a Hopf algebra over $A$ with $1$ whose counit 
$\OO_A@>>>A$ is given by $a^*\m1$ if $a=1_\l$ for some $\l\in X$ and $a^*\m0$ if $a\in\dBB$ is not of the form 
$1_\l,\l\in X$. From the definitions it is clear that $\OO_A=A\ot_\ca\OO_\ca$ as a Hopf algebra. Define an 
$A$-linear involution $\o:\OO_A@>>>\OO_A$ by $\o(a^*)=e_aa^{!*}$ (see 1.5) for $a\in\dBB$. This is an isomorphism
of the algebra $\OO_A$ onto the algebra $\OO_A$ with the opposite multiplication (see 1.5(viii)) preserving $1$; 
moreover it is compatible with the comultiplication of $\OO_A$ (see 1.5(vii)).

\subhead 3.2\endsubhead
We can reformulate 1.15(a) as follows.

(a) {\it Let $c\in\dBB$. Assume that $\l,\l',\l_1,\l_2,\l'_1,\l'_2$ in $X^+$ are such that $\l=\l_1+\l_2$, 
$\l'=\l'_1+\l'_2$, $c\in\dBB_{\l,\l'}$. There exists a function $h:\dBB_{\l_1,\l'_1}\T\dBB_{\l_2,\l'_2}@>>>\ca$ 
such that $\sum_{a\in\dBB_{\l_1,\l'_1},b\in\dBB_{\l_2,\l'_2}}h(a,b)a^*b^*=c^*$ in 
$\OO_\ca$.}
\nl
Here $a^*b^*$ is product in $\OO_\ca$.

\proclaim{Proposition 3.3} The $A$-algebra $\OO_A$ is finitely generated.
\endproclaim
We choose $\l^1,\l^2,\do,\l^r$ in $X^+-\{0\}$ such that $X^+=\NN\l^1+\NN\l^2+\do+\NN\l^r$. Let 
$\G^+=\cup_{j\in[1,r]}\dBB_{\l^j,0}$, $\G^-=\cup_{j\in[1,r]}\dBB_{0,\l^j}$, $\G=\G^+\cup\G^-$. We show:

(a) {\it The (finite) set $\{b^*;b\in\G\}$ generates the $A$-algebra $\OO_A$.}
\nl
We can assume that $A=\ca$.
For any $c\in\dBB$ let $N=N_c$ be the smallest integer $\ge0$ such that there exist $(n_1,n_2,\do,n_r)\in\NN^r$, 
$(n'_1,n'_2,\do,n'_r)\in\NN^r$ with

(b) $c\in\dBB_{\l,\l'}$ where $\l=n_1\l^1+n_2\l^2+\do+n_r\l^r$, $\l'=n'_1\l^1+n'_2\l^2+\do+n'_r\l^r$ and $\sum_jn_j+\sum_jn'_j=N$.
\nl
Clearly $N_c$ is well defined. Let $\fF$ be the $\ca$-subalgebra (with $1$) of $\OO_\ca$ generated by $\G$. It is
enough to show that $c^*\in\fF$ for any $c\in\dBB$. We shall prove this by induction on $N_c$. If $N_c=0$ then 
$c\in\dBB_{0,0}$ hence $c=1_0$ and $c^*$ is the unit element of $\OO_\ca$. If $N_c=1$ then $c\in\G$ and $c^*\in\fF$. Assume now that $N_c=N\ge2$ and that
the result is known for any $c'\in\dBB$ with $N_{c'}<N$. We can find $\l,\l'$ in $X^+$ and
$(n_j)\in\NN^r,(n'_j)\in\NN^r$ so that (b) holds. We can find $\l'_1,\l'_1,\l_2,\l'_2$ in $X^+$ such that 
$\l_1+\l_2=\l$, $\l'_1+\l'_2=\l'$ and either $(\l_1,\l'_1)=(\l^j,0)$ for some $j$ or $(\l_1,\l'_1)=(0,\l^j)$ for 
some $j$. For any $a\in\dBB_{\l_1,\l'_1}$ we have $N_a\le1$; for any $b\in\dBB_{\l_2,\l'_2}$ we have $N_b\le N-1$.
We can write $c^*=\sum_{a\in\dBB_{\l_1,\l'_1},b\in\dBB_{\l_2,\l'_2}}h(a,b)a^*b^*$ as in 3.2(a). By the induction 
hypothesis, for each $(a,b)$ in the sum we have $a^*\in\fF$, $b^*\in\fF$; hence $c^*\in\fF$. This completes the
inductive proof of (a). The proposition is proved.

\subhead 3.4\endsubhead
Let $\l\in X^+$. The sets $\dBB_{\l,0},\dBB_{0,\l}$ which enter in the definition of $\G$ in 3.3 can be described
as follows:
$$\dBB_{\l,0}=\{b^+1_{-\l};b\in\BB'_\l\},\qua \dBB_{0,\l}=\{b^-1_\l;b\in\BB_\l\}.\tag a$$
We only prove the first of these equalities; the proof of the other equality is similar. For $b\in\BB'_\l$ we have
$b^+1_{-\l}(\x_{-\l}\ot\et_0)=(b^+1_{-\l}\x_{-\l})\ot\et_0\ne0$. Since $b^+1_{-\l}\in\dBB$, it follows that 
$b^+1_{-\l}\in\dBB_{\l,0}$. Thus $\{b^+1_{-\l};b\in\BB'_\l\}\sub\dBB_{\l,0}$. Since this is an inclusion of finite
sets with the same cardinal, it is an equality, as required.

We see that 

(b) {\it the number of generators of the $A$-algebra $\OO_A$ given by the proof of 3.3 is
at most $2\sum_{s=1}^r\dim\L_{\l^s}$.}

\subhead 3.5\endsubhead
For any $a\in\dBB^+$ (resp. $a\in\dBB^-$) the restriction of the linear form $a^*:\dUU_A@>>>A$ to $\dUU^+_A$ 
(resp. $\dUU^-_A$) is denoted again by $a^*$. Let $\OO^+_A$ (resp. $\OO^-_A$) be the $A$-submodule of 
$\dUU^{+\di}_A$ (resp. $\dUU^{-\di}_A$) spanned by $\{a^*;a\in\dBB^+\}$ (resp. $\{a^*;a\in\dBB^-\}$). Note that
$\{a^*;a\in\dBB^+\}$ (resp. $\{a^*;a\in\dBB^-\}$) is an $A$-basis of $\OO^+_A$ (resp. $\OO^-_A$). We define an 
$A$-algebra structure on $\OO^{\pm}_A$ by the rule $a^*b^*=\sum_{c\in\dBB^{\pm}}\hm^{a,b}_cc^*$ for any $a,b$ in 
$\dBB^{\pm}$. (The sum is well defined by 1.16.) The (surjective) $A$-linear map $\p^+:\OO_A@>>>\OO^+_A$ given by 
$a^*\m a^*$ (if $a\in\dBB^+$), $a^*\m0$ (if $a\in\dBB-\dBB^+$), respects the algebra structures. (If $c\in\dBB^+$,
$a\in\dBB$, $b\in\dBB$, $\hm^{a,b}_c\ne0$ then $a\in\dBB^+$, $b\in\dBB^+$). Similarly, the (surjective) $A$-linear
map $\p^-:\OO_A@>>>\OO^-_A$ given by $a^*\m a^*$ (if $a\in\dBB^-$), $a^*\m0$ (if $a\in\dBB-\dBB^-$), respects the
algebra structures. It follows that the algebras $\OO^+_A$, $\OO^-_A$ are associative with $1$. Define an
$A$-linear map $\d^{\pm}:\OO^{\pm}_A@>>>\OO^{\pm}_A\ot_A\OO^{\pm}_A$ by 
$$\d^{\pm}(c^*)=\sum_{(a,b)\in\dBB^{\pm}\T\dBB^{\pm}}m_{a,b}^ca^*\ot b^*$$
for any $c\in\dBB^{\pm}$. The sum is well defined by 1.8. Note that $\p^{\pm}$ is compatible with $\d,\d^{\pm}$. 
Define an $A$-linear map $S^{\pm}:\OO^{\pm}_A@>>>\OO^{\pm}_A$ by $S^{\pm}(a^*)=s_a\ua^*$ for $a\in\dBB^{\pm}$. (We
then have $\ua\in\dBB^{\pm}$.) Note that $(\OO^{\pm}_A,\d^{\pm},S^{\pm})$ is a Hopf algebra over $A$ with $1$ and
with counit $a^*\m1$ if $a=1_\l$ for some $\l\in X$ and $a^*\m0$ if $a\in\dBB^{\pm}$ is not of the form 
$1_\l,\l\in X$. From the definitions we have $\OO^{\pm}_A=A\ot_\ca\OO^{\pm}_\ca$ as Hopf algebras. Define an 
$A$-linear isomorphism $\o:\OO^+_A@>\si>>\OO^-_A$ by $\o(a^*)=a^{!*}$ (see 1.5) for $a\in\dBB^+$. (Recall that 
$e_a=1$.) This is an isomorphism of the algebra $\OO^+_A$ onto the algebra $\OO^-_A$ with opposed multiplication.

\subhead 3.6\endsubhead
Let $\io:\OO_A@>>>\OO^-_A\ot_A\OO^+_A$ be the $A$-linear map given by the composition 
$\OO_A@>\d>>\OO_A\ot_A\OO_A@>\p^-\ot_A\p^+>>\OO^-_A\ot_A\OO^+_A$. (Note that $\io$ is an $A$-algebra 
homomorphism.)

\proclaim{Lemma 3.7} The $A$-linear map $\io$ is a split injection.
\endproclaim
For $\l\in X^+$ let $Z_\l$ be the set of pairs $(a,b)\in\dBB^-\T\dBB^+$ such
that $a\in\dBB[\l],b\in\dBB[\l],a1_\l=a,1_\l b=b$. By \cite{\LU, 4.4} for any $(a,b)\in Z_\l$ there is a unique 
$c\in\dBB[\l]$ such that $m_{a,b}^c=\d_{c,c'}$ for any $c'\in\dBB[\l]$; moreover $(a,b)\m c$ is a bijection
$\ps_\l:Z_\l@>\si>>\dBB[\l]$. Let $Z=\sqc_{\l\in X^+}Z_\l$. Then $\ps:=\sqc_\l\ps_\l:Z@>>>\dBB$ is a bijection.

Let $\uZ$ be the $A$-submodule of $\OO^-_A\ot_A\OO^+_A$ spanned by $\{a^*\ot b^*;(a,b)\in Z\}$. Define an
$A$-linear map $\r:\OO^-_A\ot_A\OO^+_A@>>>\uZ$ by sending a basis element $a^*\ot b^*$ to itself if $(a,b)\in Z$ 
and to $0$ otherwise. It is enough to show that $\r\io$ is an isomorphism of $A$-modules.

Let $c\in\dBB[\l]$, $\l\in X^+$. By definition we have $r\io(c^*)=\sum_{(a,b)\in Z}m_{a,b}^ca^*\ot b^*$. If 
$(a,b)\in Z_{\l'},\l'\in X^+$ and $m_{a,b}^c\ne0$ then as in the proof of 1.8 we have $\l'\le\l$. Thus 
$$\r\io(c^*)=a^*_0\ot b^*_0+\sum_{\l'\in X^+;\l'<\l}\sum_{(a,b)\in Z_{\l'}}m_{a,b}^c a^*\ot b^*$$
where $(a_0,b_0)=\ps_\l\i(c)\in Z_\l$. We identify $\uZ$ with $\OO_A$ as $A$-modules via the bijection 
$a^*\ot b^*\m c^*$ (with $c=\ps(a,b)$) between the bases of $\uZ$, $\OO_A$. Then we have
$$\r\io(c^*)=c^*+\sum_{\l'\in X^+;\l'<\l}\sum_{c'\in\dBB[\l']}x_{c',c}c'{}^*$$
where $x_{c',c}\in A$. We see that $\r\io$ is represented by a square upper triangular matrix with entries in $A$
and with $1$ on diagonal. It follows that $\r\io$ is an isomorphism. The lemma is proved.

\subhead 3.8\endsubhead
In the rest of this section we assume that $v=1$ in $A$. 

For $a,b\in\BB$ we have $ab=\sum_{c\in\BB}\mu_{a,b}^cc$ (in $\bof_A$) where $\mu_{a,b}^c\in A$ are zero for all 
but finitely many $c$. In the case where $A=\ca$ we define a homomorphism of $\ca$-algebras 
$\ur:\bof_\ca@>>>\bof_\ca\ot_\ca\bof_\ca$ as the restriction of the homomorphism $r:\bof@>>>\bof\ot\bof$ given in
\cite{\QG, 1.2.6}. For general $A$ we define a homomorphism of $A$-algebras $\ur:\bof_A@>>>\bof_A\ot_A\bof_A$ by 
appplying $A\ot_\ca?$ to the homomorphism $\ur:\bof_\ca@>>>\bof_\ca\ot_\ca\bof_\ca$. For any $a,b,c$ in $\BB$ we 
define $\hmu^{a,b}_c\in A$ by the following requirement: for any $c\in\BB$ we have 
$\ur(c)=\sum_{a,b\in\BB}\hmu^{a,b}_ca\ot b$ (in the last sum $\hmu^{a,b}_c$ is $0$ for all but finitely many 
$(a,b)$). We define an $A$-algebra homomorphism $S$ from $\bof_A$ to $\bof_A$ with the opposite multiplication by
$S(\th_i^{(m)})=(-1)^m\th_i^{(m)}$ for all $i\in I,m\in\NN$. For any $a\in\BB$ we have $S(a)=q_a\ta$ where 
$a\m\ta$ is an involution of $\BB$ and $q_a=\pm1$. Note that the quantities $\mu_{a,b}^c,\hmu^{a,b}_c,q_a$ (in 
$A$) are the images of the corresponding quantities in $\ca$ under the given homomorphism $\ca@>>>A$. Now 
$(\bof_A,\ur,S)$ is a Hopf algebra over $A$ with $1$ whose counit $\bof_A@>>>A$ is given by $a\m\d_{a,1}$ for 
$a\in\dBB$. For any $a\in\BB$ let $a^*:\bof_A@>>>A$ be the $A$-linear form given by $a^*(b)=\d_{a,b}$ for 
$b\in\BB$. Let $\oo_A$ be the $A$-submodule of $\bof_A^\di$ spanned by $\{a^*;a\in\BB\}$. Note that
$\{a^*;a\in\BB^+\}$ is an $A$-basis of $\oo_A$. We define an $A$-algebra structure on $\oo_A$ by the rule 
$a^*b^*=\sum_{c\in\BB}\hmu^{a,b}_cc^*$ for any $a,b$ in $\BB^{\pm}$. (The sum is well defined by homogeneity
reasons.) This algebra structure is associative with unit element $1^*$. Define an $A$-linear map 
(comultiplication) $\d_0:\oo_A@>>>\oo_A\ot_A\oo_A$ by 
$$\d_0(c^*)=\sum_{(a,b)\in\BB\T\BB}\mu_{a,b}^ca^*\ot b^*$$
for any $c\in\BB$. (The sum is well defined by homogeneity reasons.) Define an $A$-linear map $S:\oo_A@>>>\oo_A$ 
(antipode) by $S(a^*)=q_a\ta^*$ for $a\in\BB$. Note that $(\oo_A,\d_0,S)$ is a Hopf algebra over $A$ with $1$ and
with counit $a^*\m\d_{a,1}$. From the definitions we have $\oo_A=A\ot_\ca\oo_\ca$ as Hopf algebras. Define a 
(surjective) $A$-linear map $\p^{>0}:\OO_A@>>>\oo_A$ by $a^*\m b^*$ if $a=b^+1_\l$ for some $b\in\BB,\l\in X$ and
$a^*\m0$ if $a\in\dBB$ is not of the form above. Define a (surjective) $A$-linear map $\p^{<0}:\OO_A@>>>\oo_A$ by
$a^*\m b^*$ if $a=b^-1_\l$ for some $b\in\BB,\l\in X$ and $a^*\m0$ if $a\in\dBB$ is not of the form above. Note 
that $\p^{>0}$ and $\p^{<0}$ respect the $A$-algebra structures and the unit elements. They are also compatible with
$\d,\d_0$.

\subhead 3.9\endsubhead
Let $A[X]$ be the group algebra of $X$ with coefficients in $A$. Define a (surjective) $A$-linear map 
$\p^0:\OO_A@>>>A[X]$ by $a^*\m1_\l$ if $a=1_\l$ for some $\l\in X$ and $a^*\m0$ if $a\in\dBB$ is not of this form.
Note that $\p^0$ respects the $A$-algebra structures and the unit elements.

Let $\io':\OO_A@>>>\oo_A\ot_AA[X]\ot_A\oo_A$ be the $A$-linear map given by the composition 
$\OO_A@>(\d\ot1)\d>>\OO_A\ot_A\OO_A\ot_A\OO_A@>\p^{<0}\ot_A\p^0\ot_A\p^{>0}>>\oo_A\ot_AA[X]\ot_A\oo_A$. (Note 
that $\io'$ is an $A$-algebra homomorphism.) We have the following variant of Lemma 3.7:

\proclaim{Lemma 3.10} The $A$-linear map $\io'$ is a split injection. (Recall that $v=1$ in $A$.)
\endproclaim
Let $Z_\l,\ps_\l,Z,\ps$ be as in 3.7. Let $\uZ'$ be the $A$-submodule of $\oo_A\ot_AA[X]\ot_A\oo_A$ spanned by 
$\{a_1^*\ot\l\ot b_1^*;(a_1^-1_\l,1_\l b_1^+)\in Z\}$. Define an $A$-linear map
$$\r':\oo_A\ot_AA[X]\ot_A\oo_A@>>>\uZ'$$
 by sending a basis element $a_1^*\ot\l\ot b_1^*$ to itself if 
$(a_1^-1_\l,1_\l b_1^+)\in Z$ and to $0$ otherwise. It is enough to show that $\r'\io'$ is an isomorphism of 
$A$-modules. Let $c\in\dBB[\l]$, $\l\in X^+$. By definition we have 
$$\r'\io'(c^*)=\sum_{\l'\in X}\sum_{(a,b)\in Z_{\l'}}m_{a,b}^ca_1^*\ot\l'\ot b_1^*$$
where $a_1,b_1\in\BB$ are defined in terms of $a,b$ by $a=a_1^-1_\l$, $b=1_\l b_1^+$. 
If $(a,b)\in Z_{\l'}$, $\l'\in X^+$ and $m_{a,b}^c\ne0$ then as in the proof of 1.8 we have $\l'\le\l$. Thus 
$$\r'\io'(c^*)
=a^*_0\ot\l\ot b^*_0+\sum_{\l'\in X^+;\l'<\l}\sum_{(a,b)\in Z_{\l'}}m_{a,b}^ca_1^*\ot\l'
\ot b_1^*$$
where $a_0,b_0\in\BB$ are given by $(a_0^-1_\l,1_\l b_0^+)=\ps_\l\i(c)\in Z_\l$ and $a_1,b_1$ are as above. We 
identify $\uZ'$ with $\OO_A$ as $A$-modules via the bijection $a_1^*\ot\l\ot b_1^*\m c^*$ (with 
$c=\ps(a_1^-1_\l,1_\l b_1^+)$) between the bases of $\uZ'$, $\OO_A$. Then we have 
$$\r'\io'(c^*)=c^*+\sum_{\l'\in X^+;\l'<\l}\sum_{c'\in\dBB[\l']}x'_{c',c}c'{}^*$$
where $x'_{c',c}\in A$. We see that $\r'\io'$ is represented by a square upper triangular matrix with entries in 
$A$ and with $1$ on diagonal. It follows that $\r'\io'$ is an isomorphism. The lemma is proved.

\subhead 3.11\endsubhead
Let $\{x_\boc;\boc\in\NN^n\}$ be the $A$-basis of $\bof_A$ defined in 2.8. For any 
$\boc\in\NN^n$ we define an 
$A$-linear function $\x_\boc:\bof_A@>>>A$ by $\x_\boc(x_{\boc'})=\d_{\boc,\boc'}$ for any $\boc'\in\NN^n\T X$. 

\proclaim{Lemma 3.12} (a) $\{\x_\boc;\boc\in\NN^n\}$ is an $A$-basis of $\oo_A$.

(b) If $\boc',\boc''\in\NN^n$, then $\x_{\boc'}\x_{\boc''}=\x_{\boc'+\boc''}$, product in $\oo_A$.
\endproclaim
Note that each summand $\bof_{\nu,A}$ of $\bof_A$ is a free $A$-module of finite rank with basis given by its 
intersection with $\BB$; this $A$-module also has as a basis its intersection with the $A$-basis 
$\{x_\boc;\boc\in\NN^n\}$ of $\bof_A$. It follows that $\{\x_\boc;\boc\in\NN^n\}$ and $\{b^*;b\in\BB\}$ span the 
same $A$-submodule of $\bof_A^\di$ namely the set of $A$-linear functions $\bof_A@>>>A$ which vanish on 
$\bof_{\nu,A}$ for all but finitely many $\nu$. This proves (a).

We prove (b). For any $\boc\in\NN^n$ we can write uniquely $x_\boc=\sum_{b\in\BB}h_{\boc,b}b$ where $h_{\boc,b}\in A$
are zero for all but finitely many $b$; for any $b\in\BB$ we can write uniquely
$b=\sum_{\boc\in\NN^n}\tih_{b,\boc}x_\boc$ where $\tih_{b,\boc}\in A$ are zero for all but finitely many $\boc$. We 
then have $\x_\boc=\sum_{b\in\BB}\tih_{b,\boc}b^*$ and $b^*=\sum_{\boc\in\NN^n}h_{\boc,b}\x_\boc$. We have
$$\x_{\boc'}\x_{\boc''}=(\sum_{b'\in\BB}\tih_{b',\boc'}b'{}^*)(\sum_{b''\in\BB}\tih_{b'',\boc''}b''{}^*)=
\sum_{a,b',b''\in\BB}\tih_{b',\boc'}\tih_{b'',\boc''}\hmu^{b',b''}_aa^*$$
hence (b) is equivalent to the set of equalities
$$\sum_{b',b''\in\BB}\tih_{b',\boc'}\tih_{b'',\boc''}\hmu^{b',b''}_a=\tih_{a,\boc'+\boc''}.\tag c$$
for any $\boc',\boc'',a$. The equality in 2.8(a) can be rewritten as
$$\ur(x_\boc)=\sum_{\boc',\boc''\in\NN;\boc'+\boc''=\boc}x_{\boc'}\ot x_{\boc''}$$
hence as
$$\sum_{a',a'',b\in\BB}h_{\boc,b}\hmu^{a',a''}_ba'\ot a''=
\sum_{a',a''\in\BB,\boc',\boc''\in\NN^n;\boc'+\boc''=\boc}h_{\boc',a'}h_{\boc'',a''}a'\ot a''.$$
Thus.
$$\sum_{b\in\BB}h_{\boc,b}\hmu^{a',a''}_b=\sum_{\boc',\boc'')\in\NN^n;\boc'+\boc''=\boc}h_{\boc',a'}h_{\boc'',a''}$$
for any $a',a''\in\BB,\boc\in\NN^n$. We multiply both sides by 
$\tih_{c,\boc}\tih_{a',\boc'_1}\tih_{a'',\boc''_1}$ and sum over all $\boc,a',a''$. We obtain the equalities (c). This
proves (b). The lemma is proved.

\proclaim{Lemma 3.13} Assume that $v=1$ in $A$. Recall that $n=l(w_0)$. Let $A[t_1,t_2,\do,t_n]$ be the algebra of
polynomials with coefficients in $A$ in the indeterminates 
$t_1,t_2,\do,t_n$. Define an $A$-linear map 
$\k:\oo_A@>>>A[t_1,t_2,\do,t_n]$ by $\x_\boc\m t_1^{c_1}t_2^{c_2}\do t_n^{c_n}$ for any 
$\boc=(c_1,c_2,\do,c_n)\in\NN^n$. Then $\k$ is an $A$-algebra isomorphism.
\endproclaim
$\k$ is an $A$-linear isomorphism by 3.12(a). It is compatible with the algebra structures by 3.12(b). The lemma 
is proved.

\subhead 3.14\endsubhead
Define $\ti\io':\OO_A@>>>A[t_1,t_2,\do,t_n]\ot_AA[X]\ot_AA[t_1,t_2,\do,t_n]$ as the composition 
$(\k\ot1\ot\k)\io'$ where $\io':\OO_A@>>>\oo_A\ot_AA[X]\ot_A\oo_A$ is as in 3.9.

\proclaim{Theorem 3.15} Recall the assumption that $v=1$ in $A$. The $A$-linear map $\ti\io'$ is a split 
injection. It is also an imbedding of $A$-algebras with $1$. Thus, the $A$-algebra $\OO_A$ is isomorphic to a 
finitely generated $A$-subalgebra (with $1$) of $A[t_1,t_2,\do,t_n]\ot_AA[X]\ot_AA[t_1,t_2,\do,t_n]$. In 
particular, if $A$ is an integral domain then $\OO_A$ is an integral domain.
\endproclaim
The fact that $\ti\io'$ is a split injection follows from 3.10 and the fact that $\k\ot1\ot\k$ is an isomorphism 
of $A$-modules. The fact that $\ti\io'$ is an algebra homomorphism follows from the fact that $\io'$ is an algebra
homomorphism and that $\k$ is an algebra isomorphism. The fact that the algebra $\OO_A$ is finitely 
generated follows from 3.3. If $A$ is an integral domain then $A[t_1,t_2,\do,t_n]\ot_AA[X]\ot_AA[t_1,t_2,\do,t_n]$
is an integral domain hence so is $\OO_A$. The theorem is proved.

\head 4. The group $G_A$\endhead
\subhead 4.1\endsubhead 
In this section we assume that $v=1$ in $A$. 

In analogy with \cite{\KO}, we define $G_A$ as the set of $A$-algebra homomorphisms $\OO_A@>>>A$ which take $1$ to $1$. Since
$\OO_A$ is a Hopf algebra over $A$ with $1$ and with counit, $G_A$ has a natural group structure. More
explicitly, $G_A$ is the set of all $A$-linear functions $\ph:\OO_A@>>>A$ such that $\ph(a^*b^*)=\ph(a^*)\ph(b^*)$
for all $a,b\in\dBB$ and such that $\ph(1_0^*)=1$. By the correspondence $\ph\m\sum_{a\in\dBB}\ph(a^*)a$ we 
identify $G_A$ with the subset of $\hUU_A$ defined as
$$\{\sum_{a\in\dBB}n_aa\in\hUU_A;n_{1_0}=1, \sum_{c\in\dBB}\hm^{a,b}_cn_c=n_an_b\text{ for all }a,b\in\dBB\}$$
or equivalently as $\{\x\in\hUU_A;\hD(\x)=\x\ot\x,\e(\x)=1\}$. Note that $G_A$ is closed under the multiplication
in the algebra $\hUU_A$. Thus the multiplication in $\hUU_A$ induces an associative monoid structure
on $G_A$. Since the unit element $1$ of $\hUU_A$ satisfies $\hD(1)=1\ot1$, $\e(1)=1$, it belongs to $G_A$ and 
plays the role of a unit element for the monoid structure of $G_A$. If $\x\in G_A$, we have 
$\hD(S(\x))=S^{(2)}(\hD(\x))=S^{(2)}(\x\ot\x)=S(\x)\ot S(\x)$ and $\e(S(\x))=\e(\x)=1$. Thus $S(\x)\in G_A$. From
1.5(xi) we have $\mm(S^{(0,1)}(\hD(\x)))=\x$. In our case this can be written as $\mm(S^{(0,1)}(\x\ot\x))=\x$ that
is, $\mm(\x\ot S(\x))=\x$ and $\x S(\x)=\x$. We see that any element of $G_A$ has a right inverse in the monoid 
$G_A$. Hence any element of $G_A$ has a left inverse in the monoid $G_A$. Thus $G_A$ is a group (a subgroup of the
group of invertible elements of the algebra $\hUU_A$). From the definitions we see that the algebra involution 
$\o:\hUU_A@>>>\hUU_A$ restricts to an involution $\o:G_A@>>>G_A$ which is a group isomorphism.

\subhead 4.2\endsubhead 
Let $\hUU_A^{>0}$ (resp. $\hUU_A^{<0}$) be the $A$-submodule of $\hUU_A$ consisting of all elements of the form 
$\sum_{b\in\BB,\l\in X}n_b(1_\l b^+)$ (resp. $\sum_{b\in\BB,\l\in X}n_b(b^-1_\l)$) where $n_b\in A$. 
Let $\fG_A$ be the set of $A$-algebra homomorphisms $\oo_A@>>>A$ which take $1$ to $1$. Since $\oo_A$ is a 
Hopf algebra over $A$ with $1$ and with counit, $\fG_A$ has a natural group structure. Note that 
$\fG_A$ is the set of $A$-linear maps $\ph:\oo_A@>>>A$ such that $\ph(a^*b^*)=\ph(a^*)\ph(b^*)$ for all 
$a,b\in\BB$ and such that $\ph(1)=1$. (Here $a^*b^*$ is the product in $\oo_A$.) By the correspondence 
$\ph\m\sum_{a\in\BB,\l\in X}\ph(a^*)a^-1_\l$ we identify $\fG_A$ with the set of all $\x\in\hUU^{<0}_A$ such that
$\hD(\x)=\x\ot\x$ (in $\hUU^{(2)}_A$) and $\e(\x)=1$ or equivalently with $G_A^{<0}:=G_A\cap\hUU^{<0}_A$. Similarly, by the correspondence $\ph\m\sum_{a\in\BB,\l\in X}\ph(a^*)1_\l a^+$ we identify $\fG_A$ with the set of all 
$\x\in\hUU^{>0}_A$ such that $\hD(\x)=\x\ot\x$ (in $\hUU^{(2)}_A$) and $\e(\x)=1$ or equivalently with 
$G_A^{>0}:=G_A\cap\hUU^{>0}_A$. The group structure on $\fG_A$ coincides with that induced from the group structure of $G_A$ via either one of the isomorphisms $\fG_A@>\si>>G_A^{<0}$, $\fG_A@>\si>>G_A^{>0}$.

Let $\hUU_A^0$ be the set of elements of $\hUU_A$ of the form $\sum_{\l\in X}n_\l1_\l$ with $n_\l\in A$. Let
$T_A=\hUU_A^0\cap G_A$. Thus $T_A$ is the set of elements of $\hUU_A$ of the form $\sum_{\l\in X}n_\l1_\l$ with 
$n_\l\in A^\cir$ such that $\l\m n_\l$ is a group homomorphism $X@>>>A^\cir$. Clearly, $T_A$ is an abelian 
subgroup of $G_A$. Note that $T_A$ can be identified with $\Hom(X,A^\cir)=A^\cir\ot_\ZZ Y$. The isomorphism
between $A^\cir\ot Y$ and $T_A$ associates to $d\ot y$ (where $d\in A^\cir$, $y\in Y$) the element 
$\sum_{\l\in X}d^{\la y,\l\ra}1_\l$ of $T_A$. (Since $d\in A^\cir$, it can be raised to any integer power.) Note
that $T_A$ can be also identified with the set of homomorphisms of $A$-algebras $A[X]@>>>A$ preserving $1$, by 
$[g=\sum_{\l\in X}n_\l1_\l\in T_A]\m[\l\m n_\l]$. We show:

(a) {\it multiplication in $G_A$ defines an injective map $G^{<0}_A\T T_A\T G^{>0}_A@>>>A$.}
\nl
Let $\hUU_A^-$ be the set of elements of $\hUU_A$ of the form $\sum_{a\in\dBB^-}n_aa$.
Let $\x_1,\x'_1\in G^{<0}_A$, $\x_2,\x'_2\in T_A$, $\x_3,\x'_3\in G^{>0}_A$ be such that 
$\x_1\x_2\x_3=\x'_1\x'_2\x'_3$ in $G_A$. Then $\x'_3\x_3\i=\x'_2{}\i\x'_1{}\i\x_1\x_2$. The right hand side is 
contained in $\hUU^-_A$ since $G^{<0}_A,T_A$ are contained in $\hUU^-_A$ and $\hUU_A^-$ is closed under 
multiplication in $\hUU_A$. Thus we have $\x'_3\x_3\i\in\hUU^-_A$. Since $G^{>0}_A\in\hUU^{>0}_A$ and 
$\hUU^{>0}_A$ is closed under multiplication in $\hUU_A$ we have also $\x'_3\x_3\i\in\hUU^{>0}_A$. Thus 
$\x'_3\x_3\i\in\hUU^-_A\cap\hUU^{>0}_A$. From the definitions we have $\hUU^-_A\cap\hUU^{>0}_A=\{1\}$. Thus
$\x'_3\x_3\i=1$ so that $\x'_3=\x_3$. It follows that $\x_1\x_2=\x'_1\x'_2$. Then $x'_1{}\i x_1=x'_2x_2\i$ belongs
to both $\hUU_A^{<0}$ and to $\hUU_A^0$ which have intersection $\{1\}$. Thus $x'_1{}\i x_1=x'_2x_2\i=1$ and (a)
follows.

From the definitions we see that the involution $\o:G_A@>>>G_A$ interchanges $G_A^{>0}$ with $G_A^{<0}$ and its 
restriction to $T_A$ is given by $g\m g\i$.

\subhead 4.3\endsubhead 
For any $i\in I$, $h\in A$ we have $x_i(h)\in G_A^{>0}$, $y_i(h)\in G_A^{<0}$. (We use 1.18(a).)

For any $w\in W$, the elements $w',w''$ of $\hUU_A$ (see 2.5, 2.6) belong to $G_A$. (We use 2.6(a).) We have
$w'(w\i)''=1$, see 2.5. 

Now $W$ acts on $T_A$ by $w:\sum_{\l\in X}n_\l1_\l\m\sum_{\l\in X}n_{w\i(\l)}1_\l$. This corresponds to the
$W$-action on $A^\cir\ot Y$ given by $w:d\ot y\m d\ot w(y)$. This $W$ action on $T_A$ is denoted by $w:g\m w(g)$.
We show:

(a) {\it If $w\in W$, $g\in T_A$ then $w(g)=w'gw'{}\i=w''gw''{}\i$.}  
\nl
Using the definition of $w'$ we see that to prove the first equality in (a) we may assume that $w=s_i$ for some 
$i\in I$. In this case it is enough to show that $1_{s_i(\l)}=s'_i1_\l(s'_i)\i$ or that
$1_{s_i(\l)}=T'_{i,1}(1_\l)$ for any $\l\in X$. This follows from \cite{\QG, 41.1.2}. The proof of the second
equality in (a) is similar.

\subhead 4.4\endsubhead 
Let $i\in I$, $u\in A^\cir$.  Let $t_i(u)=\sum_{\l\in X}u^{\la i,\l\ra}1_\l\in T_A$. We show:

(a) $s''_i=s'_it_i(-1)$ in $G_A$.
\nl
We have
$$s'_it_i(-1)=\sum_{y\in\ZZ,x\in\ZZ,\l,\l'\in X;\l'=\l-yi',\la i,\l\ra=x+y}
(-1)^y(-1)^{\la i,\l'\ra}\th_i^{(x)-}1_\l\th_i^{(y)+}.$$
It is enough to show that $(-1)^{y+\la i,\l-yi'\ra}=(-1)^x$ in the sum above. This follows from $\la i,i'\ra=2$ 
and $y+\la i,\l\ra=x\mod2$.

\subhead 4.5\endsubhead 
Let $i\in I$. Let $e,f,g\in A$ be such that $d:=eg+f^2\in A^\cir,f\in A^\cir$. Define $e',f',g'\in A$ by
$e'=ed\i,f'=fd\i,g'=gd\i$. Note that $d':=e'g'+f'{}^2=d\i\in A^\cir$, $f'\in A^\cir$ and 
$e=e'd'{}\i,f=f'd'{}\i,g=g'd'{}\i$. We show:
$$y_i(e)t_i(f\i)x_i(g)=x_i(g')t_i(f')y_i(e').\tag a$$
We compute the left hand side in $\hUU_A$:
$$\align&\sum\Sb c,c'\in\NN,\\ \l\in X\eSb 
e^cf^{-\la i,\l\ra}g^{c'}\th_i^{(c)-}1_\l\th_i^{(c')+}=
\sum\Sb c,c'\in\NN,\l\in X;\\ \la i,\l\ra\ge c+c'\eSb
e^cf^{-\la i,\l\ra}g^{c'}\th_i^{(c)-}1_\l\th_i^{(c')+}\\&
+\sum\Sb c,c',t\in\NN,\l\in X;\\\la i,\l\ra<c+c'\eSb
e^cf^{-\la i,\l\ra}g^{c'}\bin{c+c'-\la i,\l\ra}{t}
\th_i^{(c'-t)+}1_{\l-(c+c'-t)i'}\th_i^{(c'-t)-}\\&
=\sum_{c,c'\in\NN,\l\in X;\la i,\l\ra\ge c+c'}e^cf^{-\la i,\l\ra}g^{c'}\th_i^{(c)-}1_\l\th_i^{(c')+}\\&
+\sum\Sb r,r'\in\NN,\l'\in X\\ \la i,-\l'\ra>r+r'\eSb
\sum_{t\in\NN}e^{r'+t}f^{-\la i,\l'+(r+r'+t)i'\ra}g^{r+t}\bin{-r-r'-\la i,\l'\ra}{t}
\th_i^{(r)+}1_{\l'}\th_i^{(r')-}\\&
=\sum_{c,c'\in\NN,\l\in X;\la i,\l\ra\ge c+c'}e^cf^{-\la i,\l\ra}g^{c'}\th_i^{(c)-}1_\l\th_i^{(c')+}\\&
+\sum\Sb r,r'\in\NN,\l'\in X;\\ \la i,-\l'\ra>r+r'\eSb
e^{r'}f^{-\la i,\l'+(r+r')i'\ra}g^r(1+ef^{-2}g)^{-r-r'-\la i,\l'\ra}\th_i^{(r)+}1_{\l'}\th_i^{(r')-}.\endalign$$
We compute the right hand side of (a) in $\hUU_A$:
$$\align&\sum\Sb c,c'\in\NN,\\ \l\in X\eSb
g'{}^cf'{}^{\la i,\l\ra}e'{}^{c'}\th_i^{(c)+}1_\l\th_i^{(c')-}
=\sum\Sb c,c'\in\NN,\l\in X;\\ -\la i,\l\ra>c+c'\eSb
g'{}^cf'{}^{\la i,\l\ra}e'{}^{c'}\th_i^{(c)+}1_\l\th_i^{(c')-}\\&
+\sum\Sb c,c',t\in\NN,\l\in X;\\-\la i,\l\ra\le c+c'\eSb
g'{}^cf'{}^{\la i,\l\ra}e'{}^{c'}\bin{c+c'+\la i,\l\ra}{t}\th_i^{(c'-t)-}1_{\l+(c+c'-t)i'}\th_i^{(c-t)+}\\&
=\sum_{c,c'\in\NN,\l\in X; -\la i,\l\ra>c+c'}g'{}^cf'{}^{\la i,\l\ra}e'{}^{c'}
\th_i^{(c)+}1_\l\th_i^{(c')-}\\&
+\sum\Sb s,s'\in\NN,\l'\in X;\\ \la i,\l'\ra\ge s+s'\eSb\sum_{t\in\NN}
g'{}^{s'+t}f'{}^{\la i,\l'-(s+s+t)i'\ra}e'{}^{s+t}\bin{-s-s'+\la i,\l'\ra}{t}\th_i^{(s)-}1_{\l'}\th_i^{(s')+}\\&
=\sum_{c,c'\in\NN,\l\in X; -\la i,\l\ra>c+c'}g'{}^cf'{}^{\la i,\l\ra}e'{}^{c'}\th_i^{(c)+}1_\l\th_i^{(c')-}
\\&+\sum\Sb s,s'\in\NN,\l'\in X;\\ \la i,\l'\ra\ge s+s'\eSb
g'{}^{s'}f'{}^{\la i,\l'-(s+s)i'\ra}e'{}^s(1+g'f'{}^{-2}e')^{-s-s'+\la i,\l'\ra}
\th_i^{(s)-}1_{\l'}\th_i^{(s')+}.\endalign$$
It remains to show:
$$e^cf^{-\la i,\l\ra}g^{c'}=g'{}^{c'}f'{}^{\la i,\l-(c+c')i'\ra}e'{}^{c}(1+g'f'{}^{-2}e')^{-c-c'+\la i,\l\ra},$$
$$g'{}^cf'{}^{\la i,\l\ra}e'{}^{c'}=e^{c'}f^{-\la i,\l+(c+c')i'\ra}g^c(1+ef^{-2}g)^{-c-c'-\la i,\l\ra}$$
for any $c,c'\in\NN$, $\l\in X$. This is immediate.

\subhead 4.6\endsubhead 
Let $i\in I$ and let $u\in A^\cir$. We have
$$x_i(u-1)y_i(1)x_i(u\i-1)y_i(-u)=t_i(u)\in T_A.\tag a$$
Taking $e=1$, $f=1$, $g=u\i-1$ in 4.5(a) we obtain 

$y_i(1)x_i(u\i-1)=x_i(-u+1)t_uy_i(u)$;
\nl
(a) follows.

\subhead 4.7\endsubhead 
Let $w\in W$, $i\in I$ be such that $l(ws_i)=l(w)+1$. Let $h\in A$. We show:

(a) $w''x_i(h)w''{}\i\in G_A^{>0}$.
\nl
We can find a sequence $i_1,i_2,\do,i_n$ in $I$ as in 2.8 such that $w=s_{i_1}s_{i_2}\do s_{i_{k-1}}$, $i=i_k$, 
$k-1=l(w)$. The equality 
$$T''_{i_1,1}T''_{i_2,1}\do T''_{i_{k-1},1}(\th_{i_k}^{(c)+}1_{s_{i_{k-1}}s_{i_{k-2}}\do s_{i_1}\l})=x_{c,k}^+1_\l
$$
in 2.8 (with $c\in\NN,\l\in X$, and with $x_{c,k}\in\bof_A$ independent of $\l$) can be rewritten as
$$w''(\th_{i_k}^{(c)+}1_{s_{i_{k-1}}s_{i_{k-2}}\do s_{i_1}\l})w''{}\i=x_{c,k}^+1_\l.$$
(See 2.3(c).) We multiply this by $h^c$ and sum over all $\l\in X$ and $c\in\NN$; we obtain 
$$w''x_i(h)w''{}\i=\sum_{c\in\NN}h^cx_{c,k}^+.$$
The sum in the right hand side is well defined since for some $\nu\in\NN[I]-\{0\}$ we have 
$x_{c,k}^+\in\bof_{c\nu,A}$. Thus we have $w''x_i(h)w''{}\i\in\hUU^{>0}_A$. Since $w''x_i(h)w''{}\i\in G_A$ and 
$G_A\cap\hUU^{>0}_A=G_A^{>0}$, the result follows.

\subhead 4.8\endsubhead 
In the setup of 2.8, for any $\hh=(h_1,h_2,\do,h_n)\in A^n$ we set
$$\xx_\hh=x_{i_1}(h_1)(s''_{i_1}x_{i_2}(h_2)s''_{i_1}{}\i)\do
\T(s''_{i_1}s''_{i_2}\do s''_{i_{n-1}}x_{i_n}(h_n)s''_{i_{n-1}}{}\i\do s''_{i_2}{}\i s''_{i_1}{}\i).$$
From 4.7(a) we see that $\xx_\hh\in G^{>0}_A$. We show:

(a) {\it The map $A^n@>>>G_A^{>0}$, $\hh\m\xx_\hh$ is a bijection.}
\nl
For $\hh$ as above we have
$$\align&\xx_\hh=\sum_{c_1,c_2,\do,c_n\in\NN,\l_1,\l_2,\do,\l_n\in X}h_1^{c_1}h_2^{c_2}\do h_n^{c_n}\\&
\T\th_{i_1}^{(c_1)+}1_{\l_1}T''_{i_1,1}(\th_{i_2}^{(c_2)+}1_{\l_2})
\do T''_{i_1,1}T''_{i_2,1}\do T''_{i_{n-1},1}(\th_{i_n}^{(c_n)+}1_{\l_n})\\&
=\sum_{\boc\in\NN^n,\l\in X}h_1^{c_1}h_2^{c_2}\do h_n^{c_n}x_\boc1_\l.\endalign$$
Hence for $\boc\in\NN^n$ we have
$$\xx_\hh(\x_\boc)=h_1^{c_1}h_2^{c_2}\do h_n^{c_n}.\tag b$$
Now the $A$-algebra homomorphisms $A[t_1,t_2,\do,t_n]@>>>A$ (preserving $1$) are clearly in bijection with the 
points of $A^n$; to $(h_1,h_2,\do,h_n)\in A^n$ corresponds the algebra homomorphism which takes $t_r$ to $h_r$ for
$r\in[1,n]$. Composing with $\k$ in 3.13, we obtain a bijection between the set of $A$-algebra homomorphisms 
$\oo_A@>>>A$ (preserving $1$) and the set $A^n$. From (b) we see that under this bijection $\xx_\hh\in G_A^{>0}$ 
corresponds to $\hh\in A^n$. This proves (a).

\subhead 4.9\endsubhead 
For any $\l\in X$ we define a homomorphism $\c_\l:T_A@>>>A^\cir$ by $d\ot y\m d^{\la y,\l\ra}$ for $d\in A^\cir$,
$y\in Y$ or equivalently by $\sum_{\l\in X}n_\l\l\m n_\l$. For $t\in T_A$, $h\in A$, $i\in I$, we show:
$$tx_i(h)t\i=x_i(\c_{i'}(t)h).\tag a$$
Writing $t=\sum_{\l\in X}n_\l\l$, we see that the left hand side of (a) is
$$\align&\sum_{c\in\NN,\l,\l'\in X}n_\l\i n_{\l'}h^c1_{\l'}\th_i^{(c)+}1_\l=
\sum_{c\in\NN,\l,\l'\in X}n_\l\i n_{\l'}h^c\th_i^{(c)+}1_{\l'-ci'}1_\l\\&=
\sum_{c\in\NN,\l\in X}n_\l\i n_{\l+ci'}h^t\th_i^{(c)+}1_\l=
\sum_{c\in\NN,\l\in X}n_{i'}^ch^c\th_i^{(c)+}1_\l=x_i(n_{i'}h)=x_i(\c_{i'}(t)h).\endalign$$
This proves (a).

In the setup of 2.8, let $k\in[1,n]$. Define a homomorphism $f_k:A@>>>G^{>0}_A$
by $h\m s''_{i_1}s''_{i_2}\do s''_{i_{k-1}}x_{i_k}(h)s''_{i_{k-1}}{}\i\do s''_{i_2}{}\i s''_{i_1}{}\i$,
From 3.13 we see that $f_k$ is an isomorphism of $A$ onto its image $f_k(A)$. We set 
$\l_k=s_{i_1}s_{i_2}\do s_{i_{k-1}}(i'_k)\in X$. For $t\in T_A$, $h\in A$ we show:
$$tf_k(h)t\i=f_k(\c_{\l_k}(t)h).\tag b$$
Let $t'=s_{i_{k-1}}\do s_{i_2}s_{i_1}(t)\in T_A$. Using 4.3(a), we see that (b) is equivalent to
$t'x_{i_k}(h)t'{}\i=x_{i_k}(\c_{\l_k}(t)h)$ that is, using (a), to
$x_{i_k}(\c_{i_k'}(t')h)=x_{i_k}(\c_{\l_k}(t)h)$. It is enough to show that $\c_{i_k'}(t')=\c_{\l_k}(t)$. This is
immediate from the definitions.

We show:

(c) {\it Let $w\in W$, $i\in I$. There exists $u\in A^\cir$ such that for any $h\in A$ we have
$\o(w''x_i(h)w''{}\i)=w''y_i(uh)w''{}\i$.}
\nl
Writing $w=s_{j_1}s_{j_2}\do s_{j_r}$ with $j_1,j_2,\do,j_r$ in $I$ and $r=l(w)$ we have
$$w''=s'_{j_1}t_{j_1}(-1)s'_{j_2}t_{j_2}(-1)\do s'_{j_r}t_{j_r}(-1)=w't$$
for some $t\in T_A$. Hence 
$$w''x_i(h)w''{}\i=w'tx_i(h)t\i w'{}\i=w'x_i(uh)w'{}\i$$
where $u=\c_{i'}(t)\in A^\cir$, see 4.9(a). Using 
$\o(w'')=w'$, $\o(y_i(uh))=x_i(uh))$, we obtain $w''x_i(h)w''{}\i=\o(w''y_i(uh)w''{}\i)$ and (c) is proved.

\subhead 4.10\endsubhead 
We now assume that $A$ is an algebraically closed field. Since, by 3.15, $\OO_A$ is a finitely generated 
$A$-algebra which is an integral domain, we see that $G_A$ is a connected algebraic group over $A$ with coordinate
ring $\OO_A$. Since, by 3.13,  $\oo_A$ is a finitely generated $A$-algebra which is an integral domain, we see 
that $\fG_A$ (and hence $G^{>0}_A$ and $G^{<0}_A$) are connected algebraic group over $A$ with coordinate ring 
$\oo_A$. Since the inclusions of $G^{>0}_A$ and $G^{<0}_A$ into $G_A$ correspond to surjective algebra 
homomorphisms 
$\OO_A@>>>\oo_A$ we see that $G^{>0}_A$ and $G^{<0}_A$ are closed subgroups of $G_A$. Also $T_A$ is clearly a 
subtorus of $G_A$. 

\proclaim{Theorem 4.11} In the setup of 4.10, $G_A$ is a connected reductive group over 
$A$ with associated root datum the same as the one given in 1.1.
\endproclaim
From 3.13 we see that as an affine variety, $G^{>0}_A$ is isomorphic to the affine space $A^n$.
It follows that $G^{>0}_A$ is a connected unipotent group. Similarly $G^{<0}_A$ is a connected unipotent group.

The $T_A$-action (conjugation) on the unipotent group $G^{>0}_A$ satisfies the hypotheses of 
\cite{\BO, IV,14.4}: the induced $T_A$-action on $\Lie G^{>0}$ is through one dimensional weight 
spaces corresponding to weights given by the distinct non-zero elements $\l_k (k\in[1,n])$ of $X$ (see 4.9) and
$\l_k,\l_{k'}$ are linearly independent in $X$ for $k\ne k'$. (This is a well known property of the $W$-action on
$X$.) Using {\it loc. cit} we see that any non-trivial closed subgroup of $G^{>0}_A$ which is stable under 
conjugation by $T$ contains the subgroup $f_k(A)$ (see 4.9) for some $k\in[1,n]$; moreover, the centralizer of
$T_A$ in $G^{>0}_A$ is $\{1\}$.

From 3.10 we see that the morphism $G^{<0}_A\T T_A\T G^{>0}_A@>>>G_A$ given by multiplication in $G_A$ has dense 
image. Since
$T_A$ normalizes $G^{<0}_A$ we see that $G^{>0}_AT_A$ is a closed subgroup of $G_A$. Since $G^{<0}_AT_AG^{>0}_A$ 
is the orbit of $1$ for an action of $(G^{<0}_AT)\T G^{>0}_A$ on $G_A$, we see that this orbit is locally closed 
in $G_A$. Since it is dense in $G_A$, we see that $G^{<0}_AT_AG^{>0}_A$ is open in $G_A$.

Since $G^{<0}_AT_A$ is solvable and connected, it is contained in some Borel subgroup $B$ of $G_A$. Assume that 
$G^{<0}_AT_A$ is contained properly in $B$. Now $(G^{<0}_AT_AG^{>0}_A)\cap B$ is open in $B$. If $\x_1,\x_2,\x_3$
are elements of $G^{>0}_A,T_A,G^{>0}_A$ such that $\x_1\x_2\x_3\in B$, then using $\x_1\x_2\in B$ we see that 
$\x_3\in B$. Thus $(G^{<0}_AT_AG^{>0}_A)\cap B=(G^{<0}_AT_A)(G^{>0}_A\cap B)$. If  $G^{>0}_A\cap B=\{1\}$ it would
follow that $(G^{>0}_AT_AG^{>0}_A)\cap B=G^{<0}_AT_A$ which is a proper closed subset of $B$; this contradicts the
fact that $(G^{>0}_AT_AG^{>0}_A)\cap B$ is open (non-empty) in $B$. Thus we have $G^{>0}_A\cap B\ne\{1\}$. Since 
$G^{>0}_A\cap B$ is a unipotent subgroup of $B$, it is equal to $G^{>0}_A\cap B_u$ where $B_u$ is the unipotent 
radical of $B$. Since  $G^{>0}_A\cap B_u$ is a non-trivial closed subgroup of $G^{>0}_A$, stable under conjugation
by $T_A$ we see that we have $f_k(A)\sub G^{>0}_A\cap B_u$ for some $k\in[1,n]$. Since $f_k(A)\sub G^{>0}_A$, we 
have $\o(f_k(A))\sub G^{<0}_A$. Since $G^{<0}_A$ is a unipotent subgroup of $B$, we have $G^{<0}_A\sub B_u$ hence 
$\o(f_k(A))\sub B_u$. Thus there exists $w\in W$ and $i\in I$ such that $w''x_i(A)w''\sub B_u$ and 
$\o(w''x_i(A)w'')\sub B_u$. Hence, using 4.9(c), we have $w''y_i(A)w''{}\i$. Thus $x_i(A)$ and $y_i(A)$ are 
contained in the unipotent group $w''{}\i B_uw''$. Using 4.6, we deduce that for any $u\in A^\cir$ we have
$t_i(u)\in w''{}\i B_uw''$. Since $t_i(u)\in T_A$ is semisimple it follows that $t_i(u)=1$ for all $u\in A-\{0\}$.
From the definition of $t_i(u)$ we see that if we choose $u\in A-\{0,1\}$ (recall that $A$ is algebraically 
closed) then $t_i(u)\ne1$. We have a contradiction. We see that $B=G^{<0}_AT_A$ so that $G^{<0}_AT_A$ is a Borel 
subgroup of $G_A$. Applying $\o$ we see that $\o(G^{<0}_AT_A)=G^{>0}_AT_A$ is a Borel subgroup of $G_A$. From 
4.2(a) we see that $(G^{<0}_AT_A)\cap(G^{>0}_AT_A)=T_A$. Since the unipotent radical of $G_A$ is contained in both
Borel subgroups $G^{<0}_AT_A,G^{>0}_AT_A$, it must be also contained in their intersection, the torus $T_A$. Thus,
the unipotent radical of $G_A$ is $\{1\}$. We see that $G_A$ is reductive. Since $T_A$ is the intersection of two 
Borel subgroups of $G_A$, it must be a maximal torus of $G_A$. From the definitions we see that the weights of
the adjoint action of $T_A$ on the Lie algebra of $G_A$ are exactly $\l_k (k\in[1,n])$, their negatives (each with
multiplicity $1$) and the $0$ weight with multiplicity $\dim T_A$. It follows that the root datum associated to
$G_A$ is exactly the root datum given in 1.1. The theorem is proved.

\head 5. From enveloping algebras to modified enveloping algebras\endhead
\subhead 5.1\endsubhead
In this subsection we recall some definitions from \cite{\KO}.

Let $\UU_\QQ$ be the $\QQ$-algebra with $1$ generated by the symbols $x^+,x^-$ with $x\in\bof_\QQ$ and $\uy$ with
$y\in Y$; these symbols are subject to the following relations:

$\bof_\QQ@>>>\UU_\QQ$, $x\m x^{\pm}$ is a $\QQ$-algebra homomorphism preserving $1$;

$Y@>>>\UU_\QQ$, $y\m\uy$ is $\ZZ$-linear;

$\uy\uy'=\uy'\uy$ for $y,y'\in Y$;

$\uy\th_i^{\pm}=\th_i^{\pm}(\uy\pm\la y,i\ra)$ for $y\in Y,i\in I$;

$\th_i^+\th_j^--\th_j^-\th_i^+=\d_{i,j}\ui$ for $i,j\in I$.
\nl
Thus $\UU_\QQ$ is the universal enveloping algebra attached to our root datum. There is a unique algebra 
homomorphism $\D:\UU_\QQ@>>>\UU_\QQ\ot_\QQ\UU_\QQ$ such that $\D(\uy)=\uy\ot1+1\ot\uy$ for $y\in Y$, 
$\D(\th_i^{\pm})=\th_i^{\pm}\ot1+1\ot\th_i^{\pm}$ for $i\in I$. There is a unique algebra homomorphism $S$ from 
$\UU_\QQ$ to $\UU_\QQ$ with the opposed multiplication such that $S(\uy)=-\uy$ for $y\in Y$, 
$S(\th_i^{\pm})=-\th_i^{\pm}$ for $i\in I$. There is a unique algebra homomorphism $\e:\UU_\QQ@>>>\QQ$ preserving
$1$ such that $\e(\uy)=0$ for $y\in Y$, $\e(\th_i^{\pm})=0$ for $i\in I$. Then $(\UU_\QQ,\D,S)$ is a Hopf algebra 
over $\QQ$ with $1$ and with counit $\e$.

Let $\fC'_\QQ$ be the category whose objects are unital $\UU_\QQ$-modules $M$ such that $M=\op_{\l\in X}M^\l$ (as
a vector space) where for any $\l\in X$ we have $M^\l=\{z\in M;\uy z=\la y,\l\ra z\text{ for any }y\in Y\}$.

For any $M\in\fC'_\QQ$ let $C'_M$ be the $\QQ$-subspace of $\UU_\QQ^\di$ spaned by the functions $u\m z'(uz)$ for 
various $z\in M,z'\in M^\di$. Let $\OO'_\QQ=\sum_{M\in\fC'_\QQ}C'_M$, a $\QQ$-subspace of $\UU_\QQ^\di$. Now
$\OO'_\QQ$ is a Hopf agebra over $\QQ$ in which the multiplication is the "transpose" of
$\D:\UU_\QQ@>>>\UU_\QQ\ot_\QQ\UU_\QQ$, the comultiplication is the "transpose" of the multiplication
$\UU_\QQ\ot_\QQ\UU_\QQ@>>>\UU_\QQ$, the antipode is the "transpose" of $S:\UU_\QQ@>>>\UU_\QQ$, the unit is the 
counit of $\UU_\QQ$ and the counit is given by $f\m f(1)$.

Let $\UU_\ZZ$ be the subring of $\UU_\QQ$ generated by the elements 

$\bin{\uy}{k}=\uy(\uy-1)\do(\uy-k+1))/k!$
\nl
for
$y\in Y,k\in\NN$ and $\th_i^{(m)+}$, $\th_i^{(m)-}$ for $i\in I$, $m\in\NN$. This is a $\ZZ$-lattice in 
$\UU_\QQ$, called the Kostant $\ZZ$-form of $\UU_\QQ$. Let $\OO'_\ZZ=\{f\in\OO'_\QQ;f(\UU_\ZZ)\sub\ZZ\}$. This is
a $\ZZ$-subalgebra with $1$ of $\OO'_\QQ$; it inherits from $\OO'_\QQ$ a comultiplication, an antipode and a 
counit, which make it into a Hopf algebra over $\ZZ$ (here we use the fact that $\OO'_\ZZ$ is a torsion free
$\ZZ$-module).

If $v=1$ in $A$ we set $\OO'_A=A\ot\OO'_\ZZ$; this is naturally a Hopf algebra over $A$.
The following results show that our results about $\OO_A$ apply also to $\OO'_A$.

\proclaim{Proposition 5.2} If $v=1$ in $A$, we have canonically $\OO'_A=\OO_A$ as Hopf 
algebras over $A$.
\endproclaim
The proof is given in 5.9.

\subhead 5.3\endsubhead 
Let $\PP$ be the $\QQ$-algebra of functions $X@>>>\QQ$ generated by the functions
$\uy:\l\m\la y,\l\ra$ for various $y\in Y$. We regard $\PP$ as a $\QQ$-subalgebra with $1$ of 
$\hUU_\QQ^0$ by $\ph\m\sum_{\l\in X}\ph(\l)1_\l$. 
For $y\in Y$, $h\in\ZZ$ and $k\in\NN$ we set 
$$\bin{\uy+h}{k}=(\uy+h)(\uy+h-1)\do(\uy+h-k+1))/k!\in\PP.$$
Let $\PP_\ZZ$ be the subring of $\PP$ generated by the elements $\bin{\uy}{k}$ for various $y\in Y,k\in\NN$. 
Note that

(a) $\PP_\ZZ\sub\hUU_\ZZ^0$,

(b) $\bin{\uy+h}{k}\in\PP_\ZZ$ for any $y\in Y$, $h\in\ZZ$, $k\in\NN$. 
\nl
Let $\{y_1,y_2,\do,y_r\}$ be a basis of the $\ZZ$-module $Y$. The following result is well known.

(c) {\it The elements $\ph_\kk=\bin{\uy_1}{k_1}\bin{\uy_2}{k_2}\do\bin{\uy_r}{k_r}$ (with 
$\kk=(k_1,k_2,\do,k_r)\in\NN^r$) form a $\ZZ$-basis of the $\ZZ$-module $\PP_\ZZ$ and a $\QQ$-basis of the 
$\QQ$-vector space $\PP$.}
\nl
We show:

(d) {\it Let $N\ge1$. Let $X_N=\{\l\in X;-N\le\la y_j,\l\ra\le N\text{ for }j=1,2,\do,r\}$. Let $\cm_N$ be the 
set of all elements $\sum_{\l'\in X}n_{\l'}1_{\l'}\in\hUU^0_\ZZ$ such that $n_{\l'}=0$ for $\l'\in X_N$. Let 
$\l\in X$. Then $1_\l\in\PP_\ZZ+\cm_N$.}
\nl
If $\l\n X_N$ then $1_\l\in\cm_N$ and the result holds. Now assume that $\l\in X_N$. We have
$$\prod_{j\in[1,r]}\bin{y_j+2N-\la y_j,\l\ra}{2N}=
\sum\Sb\l'\in X;\la y_j,\l'\ra\in(\iy,-\la y_j,\l\ra-2N-1]\cup
[\la y_j,\l\ra,\iy)\\
\text{ for }j\in[1,r]\eSb n_{\l'}1_{\l'}\tag e$$
where $n_{\l'}\in\ZZ$ are such that $n_\l=1$. For $\l'$ in the last sum such that 
$\la y_j,\l'\ra\in(\iy,-\la y_j,\l\ra-2N-1]$ for some $j$ we have $\la y_j,\l'\ra\le-N-1$ (since 
$-N\le\la y_j,\l\ra$) hence $\l'\n X_N$; the sum over all such $\l'$ is in $\cm_N$. On the other hand, the left 
hand side of (e) is in $\PP_\ZZ$. Thus (e) implies
$$\sum_{\l'\in X_N;\la y_j,\l'\ra\ge\la y_j,\l\ra\text{ for }j\in[1,r]}n_{\l'}1_{\l'}\in\PP_\ZZ+\cm_N$$
where $n_{\l'}\in\ZZ$ are such that $n_\l=1$. For any $\l'\in X_N$ we set 
$q_{\l'}=\sum_{j\in[1,r]}\la y_j,\l'\ra$
so that $-Nr\le q_{\l'}\le Nr$. We have
$$1_\l+\sum_{\l'\in X_N;\la y_j,\l'\ra\ge\la y_j,\l\ra\text{ for }j\in[1,r],
q_{\l'}>q_\l}n_{\l'}1_{\l'}\in
\PP_\ZZ+\cm_N.$$
This shows by descending induction on $q_\l$ (starting with $q_\l=Nr$) that $1_\l\in\PP_\ZZ+\cm_N$ for any 
$\l\in X_N$. This proves (d).

\subhead 5.4\endsubhead
There is a unique algebra homomorphism $\UU_\QQ@>>>\hUU_\QQ$ preserving $1$ such that
$\uy\m\sum_{\l\in X}\la y,\l\ra1_\l$ for any $y\in Y$, $\th_i^{\pm}\m\th_i^{\pm}$ (as in 2.7) for any $i\in I$. 
This homomorphism is injective and we identify $\UU_\QQ$ with a subalgebra of $\hUU_\QQ$ via this homomorphism. We
identify $\PP$ with a subalgebra of $\UU_\QQ$ by $\uy\m\uy$. This makes $\PP$ a subalgebra of $\hUU^0_\QQ$ (in 
the same way as in 5.3). Note that:

(a) {\it the elements $b^+\ph_\kk b'{}^-$ with $b,b'\in\BB$ and $\kk\in\NN^r$ form a $\ZZ$-basis of $\UU_\ZZ$ and
a $\QQ$-basis of $\UU_\QQ$.}

We show that

(b) $\UU_\ZZ\sub\hUU_\ZZ$.
\nl
By (a) it is enough to show that for $b\in\BB$ we have $b^+\in\hUU_\ZZ$ (this is clear),
$b^-\in\hUU_\ZZ$ 
(this is clear) and that for $\kk\in\NN$ we have $\ph_\kk\in\hUU^0_\ZZ$ (this follows from 5.3(a)). 

The category $\fC_\QQ$ (see 1.6) is equivalent to the category $\fC'_\QQ$: if $M\in\fC_\QQ$ then $M$ can be 
regarded as a $\hUU_\QQ$-module as in 1.12 and then as $\UU_\QQ$-module via the imbedding $\UU_\QQ\sub\hUU_\QQ$. 
This $\UU_\QQ$-module structure makes $M$ into an object of $\fC'_\QQ$. This gives a functor $\fC_\QQ@>>>\fC'_\QQ$
which is easily seen to be an equivalence of categories. It is well known that 

(c) {\it any object in $\fC_\QQ$ is semisimple and the simple objects of $\fC_\QQ$ are exactly the objects 
$\L_{\l,\QQ}$ for various $\l\in X^+$.}

For any $M\in\fC_\QQ$ let $C_M$ be the $\QQ$-subspace of $\dUU_\QQ^\di$ spaned by the functions $u\m z'(uz)$ for 
various $z\in M,z'\in M^\di$.

\proclaim{Lemma 5.5} Let $f\in\dUU_\QQ^\di$. The following three properties are equivalent:

(i) $f\in\OO_\QQ$;

(ii) $f\in\sum_{M\in\fC_\QQ}C_M$.

(iii) $f\in\sum_{\l\in X^+}C_{\L_{\l,\QQ}}$.
\endproclaim
Assume that (i) holds. To show (ii) we can assume that $f=a^*$ where $a\in\dBB_{\l,\l'}$ for some $\l,\l'$ in 
$X^+$. Let $M={}^\o\L_{\l,\QQ}\ot_\QQ\L_{\l',\QQ}$. Let $z=\xi_\l\ot\et_{\l'}\in M$. Define $z'\in M^\di$ by 
$z'(a'z)=\d_{a,a'}$ for any $a'\in \dBB_{\l,\l'}$. For $a''\in\dBB$ we have $z'(a''z)=a^*(a'')$ (if 
$a''\n\dBB_{\l,\l'}$ both sides are $0$ since $a''z=0$; if $a''\in\dBB_{\l,\l'}$, both sides are $\d_{a,a''}$.) 
Hence $z'(uz)=a^*(u)$ for all $u\in\dUU_\QQ$. Hence $a^*\in C_M$, as required.

Assume that (iii) holds. To show (i) we may assume that there exist $z\in M$, $z'\in M^\di$ with $M=\L_{\l,\QQ}$,
$\l\in X^+$ such that $f(u)=z'(uz)$ for all $u\in\dUU_\QQ$. Recall that, if 
$b\in\dBB-\cup_{\l';\l\ge\l'}\dBB[\l']$, then $b$ acts as $0$ on $M$ hence $f(b)=0$. Let 
$f_1=f-\sum_{b\in\cup_{\l';\l\ge\l'}\dBB[\l']}f(b)b^*\in\dUU_\QQ^\di$. Note that $f_1(b)=0$ if 
$b\in\cup_{\l';\l\ge\l'}\dBB[\l']$. If $b'\n\cup_{\l';\l\ge\l'}\dBB[\l']$ then both $f$ and 
$\sum_{b\in\cup_{\l';\l\ge\l'}\dBB[\l']}f(b)b^*$ vanish on $b'$. Hence $f_1(b')=0$. We see that $f_1=0$ and 
$f=\sum_{b\in\cup_{\l';\l\ge\l'}\dBB[\l']}f(b)b^*$. Thus $f\in\OO_\QQ$, as required.

Assume that (ii) holds. Then (iii) holds since $C_{M'\op M''}=C_{M'}+C_{M''}$ for any $M',M''$ in $\fC_\QQ$ and
any $M\in\fC_\QQ$ is a direct sum of finitely many objects of the form $\L_{\l,\QQ}$, $\l\in X^+$. The lemma is 
proved. 

\subhead 5.6\endsubhead
As in the proof of Lemma 5.5 we see that for $f\in\UU_\QQ^\di$ the following two properties are equivalent:

(i) $f\in\sum_{M\in\fC_\QQ}C'_M$.

(ii) $f\in\sum_{\l\in X^+}C'_{\L_{\l,\QQ}}$.
\nl
(Recall that $\fC_\QQ=\fC'_\QQ$.)

\subhead 5.7\endsubhead
Let $\co=\op_{\l\in X^+}\L_{\l,\QQ}\ot_\QQ\L_{\l,\QQ}^\di$, a $\QQ$-vector space. We show: 

(a) {\it the linear map $\a:\co@>>>\dUU_\QQ^\di$ whose restriction to $\L_{\l,\QQ}\ot_\QQ\L_{\l,\QQ}^\di$ is 
$z\ot z'\m[u\m z'(uz)]$ is injective.}
\nl
For any finite subset $F$ of $X^+$ let $\co_F=\op_{\l\in F}\L_{\l,\QQ}\ot_\QQ\L_{\l,\QQ}^\di$. It is enough to 
verify that for any $F$ as above, the linear map $\co_F@>>>\dUU^\di$ (restriction of $\co@>>>\dUU^\di$ in (a)) is
injective. 

Let $(\x_\l)_{\l\in F}\in\ker(\co_F@>>>\dUU^\di)$. For $\l\in F$ we have $\x_\l=\sum_jz^\l_j\ot z'_j{}^\l$ with 
$z_j^\l\in\L_{\l,\QQ}$, $z'_j{}^\l\in\L_{\l,\QQ}^\di$. For any $u\in\dUU_\QQ$ we have
$\sum_{\l\in F}\sum_jz'{}^\l_j(uz^\l_j)=0$.
Using 5.4(c) we see that for any collection $(e_\l)_{\l\in F}$, $e_\l\in\End(\L_{\l,\QQ})$ we can find 
$u\in\dUU_\QQ$ such that $u$ acts on $\L_{\l,\QQ}$ as $e_\l$. Hence we have 
$\sum_{\l\in F}\sum_jz'{}^\l_j(e_\l z^\l_j)=0$ for any $(e_\l)_{\l\in F}$ as above. Hence for any $\l\in F$ we 
have $\sum_jz'{}^\l_j(ez^\l_j)=0$  for any $e\in\End(\L_{\l,\QQ})$. Hence $\sum_jz^\l_j\ot z'{}^\l_j=0$ for any 
$\l\in F$. Thus $\x_\l=0$ for any $\l\in F$. This proves (a).

An entirely similar proof gives the following result.

(b) {\it the linear map $\a':\co@>>>\UU_\QQ^\di$ whose restriction to $\L_{\l,\QQ}\ot_\QQ\L_{\l,\QQ}^\di$ is 
$z\ot z'\m[u\m z'(uz)]$ is injective.}

From 5.5, 5.6, we see that the image of $\a$ is $\OO_\QQ$ and the
image of $\a'$ is $\OO'_\QQ$. Using (a),(b) we see that
$\a,\a'$ define isomorphisms $\co@>\si>>\OO_\QQ$, $\co@>\si>>\OO'_\QQ$. 

Now let $\hat\a:\co@>>>\hUU_\QQ^\di$ be the linear map whose restriction to $\L_{\l,\QQ}\ot_\QQ\L_{\l,\QQ}^\di$ is
$z\ot z'\m[u\m z'(uz)]$. Let $\g:\hUU^*_\QQ@>>>\dUU_\QQ^\di$, $\g':\hUU^*_\QQ@>>>\UU_\QQ^\di$ be the linear maps
transpose to the imbeddings $\dUU_\QQ@>>>\hUU_\QQ$, $\UU_\QQ@>>>\hUU_\QQ$. From the definitions we have
$\g\hat\a=\a$, $\g'\hat\a=\a'$. Since $\a$ is injective we see that $\hat\a$ is injective. Hence $\hat\a$ defines
an isomorphism $\co@>>>\hOO_\QQ$ where $\hOO_\QQ$ is the image of $\hat\a$. We see that the image of $\OO_\QQ$ 
under $\g$ is equal to the image of $\OO'_\QQ$ under $\g'$ and both these images are equal to $\hOO_\QQ$. Now
$\g,\g'$ restrict to isomorphisms $\OO_\QQ@>\si>>\hat\OO_\QQ$, $\OO'_\QQ@>\si>>\hat\OO_\QQ$. Thus both
$\OO_\QQ,\OO'_\QQ$ can be canonically identified with $\hOO_\QQ$ hence also with each other. From the definitions
we see that this identification of $\OO_\QQ,\OO'_\QQ$ is compatible with the Hopf algebra structures. 

\proclaim{Lemma 5.8} Let $F$ be a finite subset of $X^+$.

(i) Let $u\in\dUU_\ZZ$. There exists $u'\in\UU_\ZZ$ such that $u,u'$ act in the same way on 
$\L_{\l,\QQ}$ for any $\l\in F$.

(ii) Let $u'\in\UU_\ZZ$. There exists $u\in\dUU_\ZZ$ such that $u,u'$ act in the same way on 
$\L_{\l,\QQ}$ for any $\l\in F$.
\endproclaim
We can find $N\ge1$ such that $\{\l\in X;1_\l\L_{\l,\QQ}\ne0\}\sub X_N$ (notation of 5.3(d)).

To prove (i) we may assume that $u=b^+1_\l b'{}^-\in\dUU_\ZZ$ where $b,b'\in\BB$, $\l\in X$. By 5.3(d) we can 
find $u_1\in\PP_\ZZ$, $u_2\in\hUU^0_\ZZ$ such that $1_\l=u_1+u_2$ and $u_2$ acts as $0$ on $\L_{\l,\QQ}$ for any
$\l\in F$. Then $b^+u_2b'{}^-\in\hUU_\ZZ$ acts as $0$ on $\L_{\l,\QQ}$ for any $\l\in F$. Hence 
$u=b^+u_1b'{}^-+b^+u_2b'{}^-$ acts in the same way as $b^+u_1b'{}^-$ on $\L_{\l,\QQ}$ for any $\l\in F$. Hence 
$u'=b^+u_1b'{}^-\in\UU_\ZZ$ is as required in (i). 

We prove (ii). We set $u=u'\sum_{\l\in X_N}1_\l\in\dUU_\QQ$. Since $u'\in\hUU_\ZZ$ (see 
5.4(b)) we see that $u\in\hUU_\ZZ$. Thus, $u\in\dUU_\QQ\cap\hUU_\ZZ=\dUU_\ZZ$. Clearly, $u$ is as required by (ii). The lemma is proved.

\subhead 5.9\endsubhead
We prove Proposition 5.2. We can assume that $A=\ZZ$.
Let $Z$ (resp. $Z'$) be the set of all $(\x_\l)_{\l\in X^+}\in\co$, 
$\x_\l=\sum_j z^\l_j\ot z'_j{}^\l$ with $z_j^\l\in\L_{\l,\QQ}$, $z'_j{}^\l\in\L_{\l,\QQ}^\di$, such that for any 
$u\in\dUU_\ZZ$ (resp. any $u\in\UU_\ZZ$) we have $\sum_{\l\in X^+}\sum_jz'{}^\l_j(uz^\l_j)\in\ZZ$.

Under the identification $\co=\OO_\QQ$ (resp. $\co=\OO'_\QQ$) induced by $\a$ (resp. $\a'$) in 5.7, the subset
$\OO_\ZZ$ (resp. $\OO'_\ZZ$) corresponds to the subset $Z$ (resp. $Z'$) of $\co$.

Since any element $(\x_\l)_{\l\in X^+}\in\co$ has only finitely many non-zero components we see, using 5.8(i), 
that $Z'\sub Z$ and, using 5.8(ii), that $Z\sub Z'$. Thus $Z=Z'$. It follows that under the identification
$\OO_\QQ=\OO'_\QQ$, $\OO_\ZZ$ corresponds to $\OO'_\ZZ$. Proposition 5.2 follows.


\widestnumber\key{AB}
\Refs
\ref\key{\BO}\by A.Borel\book Linear algebraic groups\publ W.A.Benjamin,inc.\yr1969\publaddr New York and 
Amsterdam\endref
\ref\key{\CH}\by C.Chevalley\paper Sur certains groupes simples\jour Tohoku Math.J.\vol7\yr1955\pages14-66\endref
\ref\key{\CHH}\by C.Chevalley\paper Certains sch\'emas de groupes semi-simples\inbook S\'em. Bourbaki 1960/61\publ Soc. Math. France \yr1995 \endref
\ref\key{\JO}\by A.Joseph\book Quantum groups and their primitive ideals\bookinfo Ergebnisse der Mathematik und 
ihrer Grenzgebiete, Band 29\publ Springer Verlag\yr1995\endref
\ref\key{\KO}\by B.Kostant\paper Groups over $\ZZ$\inbook Algebraic Groups and Their Discontinuous Subgroups
\bookinfo Proc. Symp. Pure Math.\vol8 Amer. Math. Soc.\yr1966\pages90-98\endref
\ref\key{\QG}\by G.Lusztig\book Introduction to quantum groups\bookinfo Progress in Math.\vol110\publ Birkh\"auser
\yr1993\endref
\ref\key{\LU}\by G.Lusztig\paper Quantum groups at $v=\infty$\inbook Functional analysis on the eve of the 21st 
century, I\bookinfo Progr.in Math.\vol131\publ Birkh\"auser\publaddr Boston\yr1995\pages199-221\endref
\ref\key{\SO}\by Y.S.Soibelman\paper The algebra of functions on a compact quantum group and its irreducible 
representations\jour Leningrad Math.J.\vol2\yr1991\pages161-178\endref
\endRefs
\enddocument